\documentclass[%
  onecolumn
   , hidempi
]{mpi2015-cscpreprint}

\usepackage[american]{babel}

\usepackage{graphicx}
\usepackage{subcaption}
\usepackage{psfrag}

\usepackage{amssymb}
\usepackage{amsthm}

\usepackage{todonotes}
\usepackage{algorithm, algorithmic}

\input{mymymacros.sty}

\newcommand{\betat}{\ensuremath{\beta_{t}}}
\newcommand{\betaphi}{\ensuremath{\beta_{\phi}}}
\newcommand{\deltaphi}{\ensuremath{\delta_{\phi}}}

\newcommand{\Deltamu}{\ensuremath{\Delta(\p)}}

\newcommand{\DeltamustarRB}{\ensuremath{\Delta_{\text{RB}}(\p^*)}}
\newcommand{\DeltamustarEI}{\ensuremath{\Delta_{\text{EI}}(\p^*)}}

\newcommand{\nEI}{\ensuremath{n_{\text{EI}}}}

\newcommand{\nEIo}{\ensuremath{n_{\text{EI},0}}}
\newcommand{\cRB}{\ensuremath{c_{\text{RB}}}}
\newcommand{\cEI}{\ensuremath{c_{\text{EI}}}}

\newtheorem{remark}{Remark}
\begin{document}
  

\title{Fast and Reliable Reduced-Order Models for Cardiac Electrophysiology}

\author[$\ast$]{Sridhar Chellappa}
\author[$\dagger$]{Bar{\i}\c{s} Cans{\i}z}  

\author[$\ast$]{\\Lihong Feng}
\author[$\ast$]{Peter Benner}  
\author[$\dagger$]{Michael Kaliske}

\affil[$\ast$]{Max Planck Institute for Dynamics of Complex Technical Systems, 39106 Magdeburg, Germany.\authorcr
	\vspace{1mm}
	\email{chellappa@mpi-magdeburg.mpg.de}, \orcid{0000-0002-7288-3880} \authorcr
	\email{feng@mpi-magdeburg.mpg.de}, \orcid{0000-0002-1885-3269} \authorcr
	\email{benner@mpi-magdeburg.mpg.de}, \orcid{0000-0003-3362-4103} \vspace{2mm}}

\affil[$\dagger$]{Institute for Structural Analysis, Technische Universit\"at
	Dresden, 01062 Dresden, Germany.\authorcr
	\vspace{1mm}
	\email{baris.cansiz@tu-dresden.de}, \orcid{0000-0002-7568-1918}
	\authorcr
	\email{michael.kaliske@tu-dresden.de}, \orcid{0000-0002-3290-9740}}

\shorttitle{MOR for Cardiac Electrophysiology}
\shortauthor{S. Chellappa et al.}
\shortdate{}
  
\keywords{cardiac electrophysiology, model order reduction, error estimation}

  
\abstract{%
	Mathematical models of the human heart are increasingly playing a vital role in understanding the working mechanisms of the heart, both under healthy functioning and during disease. The aim is to aid medical practitioners diagnose and treat the many ailments affecting the heart. Towards this, modelling cardiac electrophysiology is crucial as the heart's electrical activity underlies the contraction mechanism and the resulting pumping action. The governing equations and the constitutive laws describing the electrical activity in the heart are coupled, nonlinear, and involve a fast moving wave front, which is generally solved by the finite element method. The simulation of this complex system as part of a virtual heart model is challenging due to the necessity of fine spatial and temporal resolution of the domain. Therefore, efficient surrogate models are needed to predict the dynamics under varying parameters and inputs. In this work, we develop an adaptive, projection-based surrogate model for cardiac electrophysiology. We introduce an a posteriori error estimator that can accurately and efficiently quantify the accuracy of the surrogate model. Using the error estimator, we systematically update our surrogate model through a greedy search of the parameter space. Furthermore, using the error estimator, the parameter search space is dynamically updated such that the most relevant  samples get chosen at every iteration. The proposed adaptive surrogate modelling technique is tested on three benchmark models to illustrate its efficiency, accuracy, and ability of generalization.
  }
 \novelty{We use an efficient and reliable a posteriori output error estimator for obtaining a reduced-order model of the cardiac electrophysiology equations. An adaptive algorithm to iteratively update the basis vectors and the parameter training set is proposed. Through this, we are able to identify reliable and fast surrogate models for different benchmark geometries.}

\maketitle

\section{Introduction}%
\label{sec:intro}
The human cardiovascular system is highly complex, exhibiting multi-scale behaviour involving multi-physics phenomena~\cite{Sachse_2005,Franzone2014}. While for several centuries the study of the human heart was largely a clinical or empirical science, in recent decades there has been significant progress in developing mathematical models which are consistent with clinical observations~\cite{aliev+panfilov96,tentusscher+noble+etal04,goktepe+kuhl09,cansiz+dal+etal15,cansiz+dal+etal17,cansiz+woodworth+etal21}. Such physics-based models seek to be complementary to empirical models which may be obtained via patient-specific imaging, test results or even population-wide studies. 

Mathematical models of the human heart are often parametrized to account for the many physical or geometrical variations possible in different patients. Such models allow for a detailed study of the effect that different parameters have on the functioning of the heart. As noted in~\cite{niederer_computational_2019}, there is an increasing trend towards using patient-specific diagnostic data within a mathematical model to enable an analytical approach to diagnosis and treatment, tailored to an individual. This philosophy encompasses the so-called \emph{cardiac digital twin} approach~\cite{Corral-Acero_EHJ2020,Peirlincketal2021,Loewe2022}. The cardiac digital twin is envisaged as a framework to not only provide a clinically-correct, physics-based surrogate of the human heart, but also to make use of available clinical data to fit several patient-specific parameters, e.g., cardiac fiber orientation, left-ventricle volume etc.

Cardiac electrophysiology, i.e., the reaction-diffusion problem in the heart tissue, is characterized by a set of ordinary differential equations (ODEs), describing excitation phenomena (i.e. reaction) and partial differential equations (PDEs) describing the wave propagation in the domain (i.e. diffusion), which are coupled in nature. The finite element method (FEM) is broadly applied for the solution of the reaction-diffusion problem in the heart which requires discretization of the space and time domains \cite{goktepe+kuhl09}. The problem of cardiac electrophysiology often requires a fine spatial and temporal discretization and, therefore, one might need to deal with a large number of equations whose solution could be time consuming and exhaustive. Therefore, developing efficient solution techniques is of great importance and is one of the essential goals in the virtual modelling of the heart electrical activity \cite{Woodworth21,woodworth+cansiz+etal22,Chak20}.

A key function of a cardiac digital twin is to address \emph{what-if} scenarios, by allowing for the repeated simulation of the model under different parameter configurations. If the digital twin has to provide results in real time, then some form of reduced-order modelling is critical. Projection-based reduced-order models (ROMs) for the human heart has been an active field of research in the past few years~\cite{Gerbeau2015,Ballarin2016,Manzoni2018,Pfaller2020,Yeetal22}. Of particular note is the recent work of Quarteroni and collaborators~\cite{Quarteroni_AN2017,QuaDM19,Fresca2020,Fresca2021} which has resulted in significant advance towards using surrogate or reduced models to accelerate the simulations of a variety of physics involved in the heart functioning. To the best of our knowledge, in these works and also in other works on MOR for cardiac models, the main methodology to obtain surrogate models/ROMs is through machine learning (ML) or through an application of the projection-based parametric proper orthogonal decomposition (POD) method~\cite{morQuaMN16}. 

Regarding ML-based approaches, while the speedup they offer is significant, a large offline computational investment is needed to realize the fast online inference. Beyond the cost, ML-based approaches are still in their development and a fair amount of investigation needs to be done to tune the various hyperparameters involved (such as network depth, number of neurons, learning rate etc.). As far as the parametric POD approach is concerned, there is a need to perform a sampling of the parameter domain to choose the samples where solution snapshots are collected. If this sampling is not done adequately, the ROM is not guaranteed to perform well on unseen parameter samples. Moreover, if a very fine sampling is carried out, this results in a high computational cost. In addition, both ML-based and POD-based approaches cannot guarantee the accuracy of the solution/quantity of interest provided by the ROM.

In this work, we propose an adaptive and greedy construction of a projection-based ROM for the reaction-diffusion problem of cardiac electrophysiology, in the spirit of the reduced basis method (RBM)~\cite{morQuaMN16,morHesRS16}. We make use of a recently proposed \emph{a posteriori} output error estimator~\cite{morCheFB19a} -- which is tailored for nonlinear dynamical systems -- to inform a greedy sampling of the parameter domain to obtain solution snapshots. Using the error estimator, we iteratively update the projection matrix to ensure it has as few basis vectors as required. Furthermore, we also make use of an adaptive parameter sampling scheme~\cite{morCheFB22} that dynamically updates the parameter training set based on the approximation error. To the best of our knowledge, a greedy construction of the (linear) projection basis of ROMs for cardiac electrophysiology has not been discussed in the literature so far. Our proposed error estimator is targeted towards output quantities of interest which are clinically relevant (e.g., the electrocardiogram (ECG)). This allows for a tailored approximation space, thus ensuring the ROMs have small dimension. Moreover, the accuracy of the ROM resulting from our proposed approach can be quantified precisely. This is crucial in making the cardiac ROM \emph{trustworthy}. We further note that the reaction-diffusion problem describing cardiac electrophysiology is formulated by the well-known phenomenological model of Aliev-Panfilov \cite{aliev+panfilov96} in the monodomain setting in the sense of \cite{goktepe+kuhl09}. 

This manuscript is organized as follows. In~\Cref{sec:mathematical_model}, we begin with an overview of the electrophysiology of the heart and briefly describe the  mathematical model. Our focus is limited to the \emph{monodomain model} in this work. Therefore, we describe in detail its spatial discretization using the FEM, followed by its temporal discretization. \Cref{sec:pmor} starts with a brief introduction to (Galerkin) projection-based reduced order models for the discretized monodomain equations and the efficient treatment of the nonlinear quantities involved therein. This section also introduces the adaptive greedy algorithm to construct the projection matrix. The proposed a posteriori error estimator and the parameter sampling strategies are discussed in detail. \Cref{sec:numerics} is devoted to the numerical implementation of the proposed method. We apply the adaptive algorithm to two benchmark examples of cardiac electrophysiology and illustrate its excellent performance. Finally, \Cref{sec:conclusion} summarizes our main contributions and charts out potential extensions.

\section{Cardiac electrophysiology}%
\label{sec:mathematical_model}
\subsection{Electrical activity in the human heart}
The electrical activity in the human heart originates in the sino-atrial nodes where pacemaker cells trigger an action potential which travels through the entire heart. Upon the excitation of the sino-atrial nodes, the electrical impulses propagate throughout the atria and are then directed to the ventricles through the atrioventricular node. The electrical signals are rapidly transmitted to the myocardium by the fast conduction system and the whole ventricles are depolarized, signifying a rise in the transmembrane potential above resting value. Afterwards, the fully depolarized ventricles go through a slow repolarization period in which the myocardium is recovered for the next excitation. During this depolarization-repolarization process, on the cellular level, several ions (e.g. Na$^-$, K$^+$, Cl$^-$, Ca$^{2+}$) are exchanged between the intracellular and extracellular media by the voltage-gated ion channels located in the cell membrane, leading to a significant alteration in the voltage or the so-called transmembrane potential. Moreover, the coupling between the electrophysiology and mechanical activity of the myocytes is linked through the intracellular Ca$^{2+}$ ion concentrations. The intake of Ca$^{2+}$ into the intracellular medium upon the myocyte depolarization triggers a series of events that results in myocardial contraction.

\subsection{Modelling approaches of cardiac electrophysiology}
Over the past half century, a wide variety of mathematical models have been proposed to model the electrical activity in the heart, starting from the pioneering work of Hodgkin and Huxley \cite{hodgkin+huxley90}. Two essential approaches exist for the description of cardiac electrophysiology: ionic and phenomenological models. The ionic models represent a sophisticated representation of the cardiomyocytes by considering the local evolution of individual ion species in line with experimental observations \cite{tentusscher+noble+etal04,wong+goktepe+etal11}. These models are useful when one needs to study the influence of a particular ion activity on cardiac electrophysiology, e.g., drug application. However, a high number of evolution equations for ion concentrations, ionic currents, and gating variables results in a demanding computational effort. On the other hand, for the investigation of wave propagation in healthy and pathological cases (e.g. arrhythmia and dyssynchrony) on the tissue or the organ level, the phenomenological models are convenient due to their ease of implementation and relatively less computational load compared to the ionic models. In this context, the FitzHugh-Nagumo model \cite{fitzhugh61,nagumo+arimoto+etal62} provides an excellent description of the excitable heart cells. The model is able to mimic the intrinsic characteristics of the transmembrane potential and lumps the influence of all ionic currents in a single slow recovery variable $r$. The FitzHugh-Nagumo model is often employed for the oscillatory cardiac cells (i.e., sino-atrial node or Purkinje fibres). On the other hand, the Aliev-Panfilov model \cite{aliev+panfilov96}, which is a modification of the FitzHugh-Nagumo model, was developed in order to model the non-oscillatory cardiac cells, i.e., myocytes, and has been extensively utilized by many researchers.

Concerning the propagation of electrical waves, the mathematical model of a single cardiac cell is furnished with a conduction term enabling the spatial description of travelling excitation waves for simulations performed in tissue or organ level. Then, the system consists of ODEs describing the ionic current and one PDE in case of a monodomain setting which is often employed if an external electrical field does not need to be applied. The numerical solution of cardiac electrophysiology problems is often handled by the FEM in the literature \cite{rogers+mcculloch94a,rogers02,sundnes+lines+etal05,goktepe+kuhl09}.

\subsection{Mathematical setting}
In this work, the electrical activity in the non-oscillatory myocardial tissue is modelled with coupled PDEs and ODEs describing the dynamics of the transmembrane potential, the recovery variable, and the wave propagation in a monodomain setting. The depolarization-repolarization behaviour is characterized by the well-known phenomenological model of Aliev-Panfilov \cite{aliev+panfilov96}. The system of coupled governing equations is written as
\begin{subequations}
\label{eq:pde-ode_monodomain}
\begin{align}
\label{eq:pde_monodomain}
\frac{\partial \Phi}{\partial t} &= \text{div} (\bD \cdot \grad \Phi) + F^{\phi} + i_{s}(t),\\[0.5em] 
\label{eq:ode_monodomain}
\betat \frac{\partial r}{\partial t} &= F^{r}, 
\end{align}
\end{subequations}
with $\Phi:=\Phi(\bz, t; \p) $ being the transmembrane potential and $r := r(\bz, t; \p)$ being the recovery variable, where $\bz \in \mathcal{B}$ is the spatial variable and $t \in [0, T]$ is the time variable. Further, $\p \in \R^{p}$ are the parameters of the system and $\bD=d_{\text{iso}}\mathbf{I}$ denotes the conductivity tensor in the domain with the isotropic conductivity parameter $d_{\text{iso}}$ and the identity matrix $\mathbf{I}$. 

For given parameters, we further endow the above system with initial conditions $\Phi(\bz, 0; \p) =: \Phi_{0}$ and $r(\bz, 0; \p) = \mathbf{0}$. $i_{s}(t)$ denotes a time-dependent external stimulus applied to the cardiac domain $\mathcal{B}$. Furthermore, the governing equations are subjected to Dirichlet and Neumann boundary conditions
\begin{equation}
\label{eq:bound-cond}
\Phi = \bar\Phi \quad \text{on}~\partial{\mathcal{B}}_{\Phi}, \quad r = \bar r \quad \text{on}~\partial{\mathcal{B}}_{r} \quad \text{and} \quad (\bD \cdot \grad \Phi)\cdot \mathbf{n} = \bar q  \quad \text{on}~ \partial{\mathcal{B}}_q.
\end{equation}
We denote by $t_{s}$ the starting time of the external stimulus and with $t_{e}$ we denote the end time of the applied stimulus. We further define $\delta t := t_{e} - t_{s}$ to be the duration of the applied stimulus. 
\begin{remark}
In contrast to the innovative finite element formulation proposed by G\"oktepe et al. \cite{goktepe+kuhl09}, where the recovery variable $r$ is treated as an internal variable at Gauss point level, $r$ in this study is discretized as an additional degree of freedom along with the transmembrane potential $\Phi$. The motivation behind this is to be able to construct efficient reduced-order models. Treating $r$ locally as an internal variable, while advantageous, still involves computational complexity scaling with the number of elements in the finite element discretization. However, by having $r$ as a global variable, it can be projected onto a reduced subspace leading to a significant reduction of the degrees of freedom involved. This will be illustrated in the numerical results. Other works developing ROMs for cardiac electrophysiology and cardiac electromechanics have also considered the recovery voltage as an additional degree of freedom~\cite{Pagani2018,Bonomi2017}.
\end{remark}
The quantities $F^{\phi}, F^{r}$ are obtained based on the Aliev-Panfilov model and correspond to the expressions:
\begin{subequations}
\label{eq:AP}
\begin{align}
\label{eq:AP_Fphi}
\frac{\betat}{\betaphi} F^{\phi} &= \frac{\partial \phi}{\partial \tau} = c \phi (\phi - \alpha)(1-\phi) - r \phi,\\[0.5em]
\label{eq:AP_Fr}
\betat F^{r} &= \frac{\partial r}{\partial \tau} = \left[\gamma + \frac{\mu_{1} r}{\mu_{2} + \phi }\right] \left[-r - c \phi (\phi - b - 1)\right].
\end{align}
\end{subequations}
In the above equations, the quantities $(\betat, \betaphi, \alpha, b, c, \gamma, t_{s}) := \p$ are parameters of the system. Of special interest to us is the parameter $\gamma$ which controls the repolarization of the cardiac muscles and the parameter $t_{s}$ that defines the initiation time of the applied input stimulus. Furthermore, the Aliev-Panfilov model makes use of the dimensionless variant of the transmembrane potential, denoted as $\phi$ and the dimensionless time $\tau$. They bear the relation with their dimensioned counterparts
\begin{subequations}
\label{eq:unit_transform}
\begin{align}
\label{eq:unit_transform_Phi}
\phi &= \frac{\Phi + \deltaphi}{\betaphi},\\[0.5em] 
\label{eq:unit_transform_t}
\tau &= \frac{t}{\betat}.
\end{align}
\end{subequations}

\begin{table}[t!]
\def\arraystretch{1.2}
\centering
\begin{tabular}{|c|r|c|}
\hline
Quantity & Value & Dimension\\
\hline
\hline
\betat 		& $12.9$ 	& [ms]\\
\betaphi 	& $100$ 	& [mV]\\
\deltaphi	& $-80$     & [mV]\\
$c$			& $8$		& [-]\\
$\alpha$	& $0.01$	& [-]\\
$b$			& $0.15$    & [-]\\
$\mu_{1}$   & $0.2$     & [-]\\
$\mu_{2}$   & $0.3$     & [-]\\
\hline		
\end{tabular}
\caption{Aliev-Panfilov monodomain equation fixed parameters}
\label{tab:AP_model_fixed_parameters}
\end{table}
\Cref{tab:AP_model_fixed_parameters} summarizes the values of the fixed parameters in \cref{eq:pde-ode_monodomain,,eq:AP,,eq:unit_transform} while $\gamma$ and $t_{s}$ are treated as free parameters. Therefore, for the remainder of this work, we consider either $\p := \gamma$ (when no external stimulus is applied) or $\p := (\gamma, t_{s})$.

\subsection{Spatial discretization of the monodomain equations}
The discretization of the system~\cref{eq:pde-ode_monodomain} is carried out using a Galerkin FEM. Recall that we denote the cardiac domain of interest as $\mathcal{B}$. First, we obtain the weak formulation of the system~\cref{eq:pde-ode_monodomain} within $\mathcal{B}$. To this end, we introduce the finite element space $X$. The weak formulation of the problem is then to find $\Phi(\bz, t; \p), r(\bz, t; \p) \in X$ such that
\begin{align}
\label{eq:weak_form_Phi}
\int\limits_{\mathcal{B}} \frac{\partial \Phi}{\partial t} \psi d \bz - \int\limits_{\mathcal{B}} ({\bD} \cdot \grad \Phi)  \grad \psi d \bz - \int\limits_{\mathcal{B}} F^{\phi} \psi d \bz - \int\limits_{\mathcal{B}} i_{s}(t) \psi d \bz - \int\limits_{\partial{\mathcal{B}}_q}\psi \bar{q} d\bz  &= 0,\\[0.5em]
\label{eq:weak_form_r}
\int\limits_{\mathcal{B}} \betat \frac{\partial r}{\partial t} \psi d \bz - \int\limits_{\mathcal{B}} F^{r} \psi d \bz &= 0
\end{align}
with $\psi \in X$ representing any of the trial functions, see \cite{cansiz+kaliske22} for a similar treatment. Next, consider a finite dimensional approximation $X_{N} \subset X$ of dimension $N$. Based on the choice of the FEM space used, let $\{\psi_{i}\}_{i=1}^{N}$ be a set of basis functions such that, the following approximation holds
\begin{align*}
\Phi(\bz, t; \p) \approx \sum\limits_{j=1}^{N} x^{\Phi}_{j}(t, \p)\, \psi_{j}(\bz)\,\, \text{and}\,\, r(\bz, t) \approx \sum\limits_{j=1}^{N} x^{r}_{j}(t, \p)\, \psi_{j}(\bz).
\end{align*}
We further define $\bx^{\Phi} := [x^{\Phi}_{1}, \ldots, x^{\Phi}_{N}]^{T} \in \R^{N}$  and $\bx^{r} := [x^{r}_{1}, \ldots, x^{r}_{N}]^{T} \in \R^{N}$ as the coefficient vectors for the above ansatz. Subsequently, Galerkin projection is applied to \cref{eq:weak_form_Phi,,eq:weak_form_r}  to compute $\bx^{\Phi}, \bx^{r}$ by solving the discretized monodomain equations
\begin{subequations}
\label{eq:disc_monodomain}
\begin{align}
\label{eq:disc_monodomain_Phi}
\bM \frac{d \bx^{\Phi}(t, \p)}{dt} &= \bS \bx^{\Phi}(t, \p) + \bM \bff^{\phi}(\bx^{\Phi}, \bx^{r}, \p) + \breve{\bB} i_{s}(t),\\[0.5em]
\label{eq:disc_monodomain_r}	
\betat \bM \frac{d \bx^{r}(t, \p)}{dt} &= \bM \bff^{r}(\bx^{\Phi}, \bx^{r}, \p),\\[0.8em]
\label{eq:disc_monodomain_init}
\bx^{\Phi}(0, \p) &= \bx^{\Phi}_{0}(\p)\,\, \text{and}\,\, \bx^{r}(0, \p) = \mathbf{0}.
\end{align}
\end{subequations}
In the above equations, the mass matrix is denoted $\bM \in \R^{N \times N}$ having entries
\begin{align*}
[\bM]_{ij} :=  \int\limits_{\mathcal{B}} \psi_{i} \psi_{j} d \bz, \qquad i, j = 1, 2, \ldots, N.
\end{align*}
Further, the stiffness matrix is denoted as $\bS \in \R^{N \times N}$ with entries
\begin{align*}
[\bS]_{ij} :=  \int\limits_{\mathcal{B}} \grad \psi_{i} \bD \grad \psi_{j} d \bz, \qquad i, j = 1, 2, \ldots, N.
\end{align*}
The input matrix is denoted $\breve{\bB} \in \R^{N}$ and its entries are given by
\begin{align*}
[\breve{\bB}]_{i} :=  \int\limits_{\mathcal{B}} \psi_{i} d \bz, \qquad i = 1, 2, \ldots, N.
\end{align*}
The terms $ \bff^{\phi}, \bff^{r} \in \R^{N}$ are discretized nodal values of the Aliev-Panfilov model variables $F^{\phi}$ and $F^{r}$, respectively.
To write the discretized nonlinear terms ($\bM \bff^{\phi}$ and $\bM \bff^{r}$), we utilize the Ionic Current Interpolation (ICI) approach~\cite{Pathmanathan2011}. Note that compared to the State Variable Interpolation (SVI) approach~\cite{Pathmanathan2011}, the ICI method is less accurate. However, ICI has far less computational cost. By defining $\bx := [\bx^{\Phi}\,\,\bx^{r}]^{T} \in \R^{2N}$, we can rewrite~\cref{eq:disc_monodomain} as below
\begin{subequations}
\label{eq:disc_monodomain_coupled}
\begin{align}
\bE \frac{d \bx(t, \p)}{dt} &= \bA\bx(t, \p) + \bM_{f} \bff(\bx, \p) + \bB i_{s}(t),\\
\bx(0, \p) &= \bx_{0}(\p),
\end{align}
\end{subequations}
where 
\begin{align*}
\bE &:=
\begin{bmatrix}
\bM & \mathbf{0}\\
\mathbf{0} & \betat \bM
\end{bmatrix} \in \R^{2N \times 2N}, \,\,\,
\bA :=
\begin{bmatrix}
\bS & \mathbf{0}\\
\mathbf{0} & \mathbf{0}
\end{bmatrix} \in \R^{2N \times 2N},\\[0.8em]
\bM_{f} &:= 
\begin{bmatrix}
\bM & \mathbf{0}\\
\mathbf{0} & \bM
\end{bmatrix} \in \R^{2N \times 2N}, \,\,\,
\bff :=
\begin{bmatrix}
\bff^{\phi}\\
\bff^{r}
\end{bmatrix} \in \R^{2N}
\text{and} \,\,\,
\bB :=
\begin{bmatrix}
\breve{\bB}\\
\mathbf{0}
\end{bmatrix} \in \R^{2N}.
\end{align*}
In the sequel, we refer to~\cref{eq:disc_monodomain_coupled} as the full-order model or FOM.
\begin{remark}
    The spatial discretization is implemented through \textsf{FEAP} (Finite Element Analysis Program)~\cite{taylor_feap}. We extract the discretized matrices ($\bM, \bS$ and $\breve{\bB}$) from \textsf{FEAP}.
\end{remark}

\subsection{Quantities of clinical interest}
In cardiology departments, ECG is one of the essential and indispensable tools in evaluating cardiac function. It is non-invasive, easy, fast, and cheap, yet delivers a great deal of information regarding the electrical activity in the heart. Basically, the electrical wave propagation, in other words, the electrical flux, is globally projected to some predefined principal directions in the heart and for each direction a flux diagram is recorded. These diagrams are used to investigate rhythm irregularities, structural changes, and infarcted zones in the heart \cite{einthoven12}.

\subsection{Time discretization of the monodomain equations}
Concerning the time discretization of the monodomain (and bidomain) equations, existing literature covers (a) fully implicit~\cite{Cansiz2022,Woodworth21}, (b) fully explicit~\cite{Puwal-Roth2007}, (c) implicit-explicit~\cite{Pathmanathan2011}, and (d) operator-splitting approaches~\cite{Zhilin1999_operatorsplit1,Krishnamoorthi2013}. In this work, we adopt a first-order implicit-explicit approach where we treat the nonlinear term explicitly and the diffusion/conduction term implicitly~\cite{AscRW95}. This approach offers a significant reduction in the computational cost, while yielding accurate results. We divide the time domain $[0, T]$ into $N_{t}$ nodes of fixed step size $\delta t$. The time-discrete version of \cref{eq:disc_monodomain_coupled} reads 
\begin{subequations}
\label{eq:fullydisc_monodomain_coupled}
\begin{align}
\bbE \bx^{k}(\p) &= \bbA \bx^{k-1}(\p) + \delta t \left( \bM_{f} \bff(\bx^{k-1}, \p) + \bB i_{s}^{k} \right),\\
\bx^{0}(\p) &= \bx_{0}(\p),
\end{align}
\end{subequations}
with $\bbE := \left(\bE - \delta t \bA \right) \in \R^{2N \times 2N}$, $\bbA := \bE \in \R^{2N \times 2N}$ and 
$\bff(\bx^{k-1}, \p) := 
\begin{bmatrix}
\bff^{\phi}(\bx^{\Phi^{k-1}}, \bx^{r^{k-1}}, \p)\\
\bff^{r}(\bx^{\Phi^{k-1}}, \bx^{r^{k-1}}, \p)
\end{bmatrix} \in \R^{2N}$.
Often, a very fine mesh size is required in order to accurately capture the propagation of the action potential wavefront. This leads to the dimension of the discretized system~\cref{eq:fullydisc_monodomain_coupled}, $N$, being large (ranging from a few thousands to several millions in highly resolved models). As a result, solving the system repeatedly for multiple parameters $\p$ incurs a large computational effort. This cost can be mitigated through the use of reduced-order modelling approaches.

\section{Projection-based model order reduction\\for the monodomain equations}%
\label{sec:pmor}
\subsection{Galerkin reduced-order model}
Linear projection-based model order reduction (PMOR) approaches project the system equations (such as \cref{eq:fullydisc_monodomain_coupled}) to a low-dimensional (linear) subspace, denoted by $\mathcal{V}_{n}$, of the true solution space $\mathcal{V}_{N}$, i.e., $\mathcal{V}_{n} \subset \mathcal{V}_{N}$. The underlying assumption is that the dimension of $\mathcal{V}_{n}$ denoted $n$ is significantly smaller than that of $\mathcal{V}_{N}$, or $n \ll N$. The model reduction machinery is then aimed at identifying a suitable basis for $\mathcal{V}_{n}$ denoted as $\bV$. Very recently, PMOR has been extended to cases where the system equations are projected to quadratic or nonlinear manifolds~\cite{morBarF22,morBarFM23}. The advantage here is that this typically leads to a smaller $n$, at the cost of increased computational complexity. We limit our focus in this work to linear PMOR.

Suppose the true solution of the coupled system~\cref{eq:fullydisc_monodomain_coupled} at a given parameter and time instance, i.e., $\bx^{\Phi^{k}}(\p), \bx^{r^{k}}(\p)$ can be approximated, respectively, in the subspaces $\mathcal{V}_{n}^{\Phi}$ and $\mathcal{V}_{n}^{r}$ as a linear combination of the basis functions, i.e.,
\begin{subequations}
\label{eq:linearmor_ansatz}
\begin{align}
\bx^{\Phi^{k}}(\p) \approx \btx^{\Phi^{k}} &:= \sum\limits_{i=1}^{n^{\Phi}} \bv^{\Phi}_{i} \hat{x}_{i}^{\Phi^{k}}(\p) = \bV^{\Phi} \bhx^{\Phi^{k}}(\p),\\[0.8em]
\bx^{r^{k}}(\p) \approx \btx^{r^{k}} &:= \sum\limits_{i=1}^{n^{r}} \bv^{r}_{i} \hat{x}_{i}^{r^{k}}(\p) = \bV^{r} \bhx^{r^{k}}(\p).
\end{align}
\end{subequations}
In the above ansatz, $\bV^{\Phi} \in \R^{N \times n^{\Phi}}$ is a basis for the projection subspace corresponding to the $\bx^{\Phi}$ variable while $\bV^{r} \in \R^{N \times n^{r}}$ is a basis for the projection subspace corresponding to the $\bx^{r}$ variable. We define the block diagonal matrix
\begin{align}
\label{eq:blockdiag_V}
\bV := 
\begin{bmatrix}
\bV^{\Phi} & \\
& \bV^{r}
\end{bmatrix} \in \R^{2N \times (n^{\Phi} + n^{r})}.
\end{align}
For a simplified notation, we let $n := (n^{\Phi} + n^{r})$ for the remainder of this work.

To obtain the reduced-order model corresponding to \cref{eq:fullydisc_monodomain_coupled}, we insert the ansatz \cref{eq:linearmor_ansatz} into \cref{eq:fullydisc_monodomain_coupled} and perform a Galerkin projection with $\bV$. The resulting ROM reads
\begin{subequations}
\label{eq:fullydisc_monodomain_coupled_rom}
\begin{align}
	\bbhE \bhx^{k}(\p) &= \bbhA \bhx^{k-1}(\p) + \delta t \left( \widehat{\bM}_{f} \btf^{k-1} + \bhB i_{s}^{k} \right),\\
	\bhx^{0}(\p) &= \bhx_{0}(\p).
\end{align}
\end{subequations}
In the above ROM, $\bbhE := \bV^{T} \bbE \bV \in \R^{n \times n}, \bbhA := \bV^{T} \bbA \bV \in \R^{n \times n}$. Further, $\widehat{\bM}_{f} := \bV^{T} \bM_{f} \in \R^{n \times 2N}$, $\btf^{k-1} :=  \bff^{k-1}(\btx^{k-1}, \p) \in \R^{2N}$ and $\bhx^{0}(\p) := \bV^{T} \bx^{0}(\p) \in \R^{n}$. Compared to solving the FOM~\cref{eq:fullydisc_monodomain_coupled}, the solution of~\cref{eq:fullydisc_monodomain_coupled_rom} is faster since the dimension of the ROM is much smaller, i.e., $n \ll 2N$.

\subsection{Treating the nonlinearities efficiently}
While the ROM~\cref{eq:fullydisc_monodomain_coupled_rom} is of smaller dimension, its computation is still inefficient as the quantity $\btf^{k-1}$ is still of dimension $2N$. There has been a wealth of methods proposed in the last decade to remedy this issue; these methods fall under the family of \emph{hyperreduction methods} with the empirical interpolation method (EIM)~\cite{morBarMNetal04}, the discrete empirical interpolation method (DEIM)~\cite{morChaS10} being the most widely used. Both EIM and DEIM use a linear approximation for the nonlinear term
\begin{align*}
\bff \approx \bff_{\text{EI}} := \bU \left(\bP^{T} \bU\right)^{-1} \bP^{T} \bff.
\end{align*} 
The matrix $\bU \in \R^{N \times n_{\text{EI}}}$ is an orthogonal matrix which is obtained from the nonlinear snapshots (in case of EIM) or from the POD basis of the nonlinear snapshots (in case of DEIM). The matrix $\bP \in \R^{N \times n_{\text{EI}}}$ is a selection matrix consisting of only $\{0, 1\}$ as its entries. The efficiency of hyperreduction comes from the fact that $\bU \left(\bP^{T} \bU\right)^{-1}$ can be precomputed which leaves only the small vector $\bP^{T} \bff \in \R^{n_{\text{EI}}}$ to be evaluated and since, for many systems, $n_{\text{EI}} \ll N$, this leads to significant speed up. 

We approximate the two nonlinear terms $\bff^{\phi}$ and $\bff^{r}$~in \cref{eq:fullydisc_monodomain_coupled} separately. We define $(\bU^{\phi}, \bP^{\phi})$ to be the hyperreduction quantities to approximate $\bff^{\phi}$ and  $(\bU^{r}, \bP^{r})$ to be the hyperreduction quantities to approximate $\bff^{r}$. Further, we define the block hyperreduction projection matrix as
\begin{align}
\label{eq:blockU}
\bU :=
\begin{bmatrix}
\bU^{\phi} & \\
& \bU^{r}
\end{bmatrix} \in \R^{2N \times (\nEI^{\phi} + \nEI^{r})}
\end{align}
and the selection matrix to be
\begin{align}
\label{eq:blockP}
\bP :=
\begin{bmatrix}
\bP^{\phi} & \\
& \bP^{r}
\end{bmatrix} \in \R^{2N \times (\nEI^{\phi} + \nEI^{r})}.
\end{align}
For easy notation, we define $\nEI := (\nEI^{\phi} + \nEI^{r})$.
The hyperreduced ROM corresponding to~\cref{eq:fullydisc_monodomain_coupled_rom} is given by
\begin{subequations}
\label{eq:fullydisc_monodomain_coupled_romEI}
\begin{align}
\bbhE \bhx^{k}(\p) &= \bbhA \bhx^{k-1}(\p) + \delta t \left( \widehat{\bM}_{f} \bff_{\text{EI}}^{k-1} + \bhB i_{s}^{k} \right),\\
\bhx^{0}(\p) &= \bhx_{0}(\p)
\end{align}
\end{subequations}
with $\bff_{\text{EI}}^{k-1} := \bU \left(\bP^{T} \bU\right)^{-1} \bP^{T} \bff(\btx^{k-1})$. In~\cref{eq:fullydisc_monodomain_coupled_romEI}, the quantity $\widehat{\bM}_{f}\,\bU \left(\bP^{T} \bU\right)^{-1} \in \R^{n \times \nEI}$ can be precomputed once and stored in memory. At each time step, the nonlinearity $\bff(\btx^{k-1})$ just needs to be evaluated at $\nEI$ indices (i.e., $\bP^{T} \bff(\btx^{k-1})$) which are specified through $\bP$.

\subsection{Computing the basis $\bV$}
In past works on MOR applied to cardiac models~\cite{Ballarin2016,Bonomi2017,YangVeneziani2017,Manzoni2018,Pagani2018,Pfaller2020,Khan_2022}, the predominant approach to compute the basis $\bV$ is the parametric proper orthogonal decomposition (POD). In the parametric POD approach, a sampling of the parameter space $\mathcal{P}$ is first done to select $n_{s}$ parameter samples $\{\p_{1}, \p_{2}, \ldots, \p_{n_{s}}\}$. The FOM~\cref{eq:fullydisc_monodomain_coupled} is then computed at all the selected parameter samples to obtain a snapshot matrix $\bX_{s}(\p_{i}) := \left[ \bx^{k}(\p_{i}) \right]_{i=1, k=0}^{i=s, k=N_{t}} \in \R^{N \times s\,N_{t} }$. The projection matrix $\bV$ is obtained by performing a singular value decomposition (SVD) of $\bX_{s}$ and setting $\bV := \bU_{\bX}(:, 1:n)$, where $\bU_{\bX}$ is the matrix containing the left singular vectors of the snapshot matrix $\bX_{s}$. While straightforward, this approach requires many FOM solutions at different parameter samples. For complex parameter domains, this can lead to a high offline computational effort. Moreover, the choice of the parameter samples is rather heuristic. Different sampling approaches, including, uniform, random, Latin hypercube sampling etc. are possible. Nevertheless, if the parameter space is not adequately sampled, the resulting ROM can be inaccurate. The choice of the subspace dimension $n$ in the POD approach is based on the decay of the singular values and the energy criterion. Typically, a ``low enough" tolerance is chosen and $n$ is set to be the smallest value of $n_{d}$ for which the following expression holds
\[
	n := \text{smallest}\,\, n_{d}\,\, \text{for which}\,\,\, \frac{\sum\limits_{j=n_{d}+1}^{n_{X}} \sigma_{j}^{2}}{\sum\limits_{j=1}^{n_{X}} \sigma_{j}^{2}} < \epsilon_{\text{SVD}}.
\]
Here, $\sigma_{1} \geq \sigma_{2} \geq \cdots \sigma_{j} \geq \cdots \sigma_{n_{X}} > 0$ are the $n_{X}$ non-zero singular values of the matrix $\bX_{s}$ and $\epsilon_{\text{SVD}}$ is a preferred tolerance. Such an approach to choose $n$ is heuristic, as it has no connection to the actual error of the dynamics or that of some desired output. For example, if there is a quantity of interest/output, then the cut-off of the singular values has no direct relation to the actual error in the output resulting from the ROM.

In this work, we utilize the adaptive POD-Greedy algorithm~\cite{morCheFB19a} to determine the projection basis $\bV$. This approach is an improvement over the standard parametric POD approach described above due to the following reasons
\begin{itemize}
\item The parameter sampling is done based on a greedy algorithm, driven by an \emph{a posteriori} output error estimator. This enables the sampling to be informed by the actual approximation quality of the output of interest.

\item The number of FOM solves is (less that or) equal to the number of iterations $n_{g}$ of the adaptive POD-Greedy algorithm; typically, $n_{g} \ll n_{s}$, therefore, this leads to a more efficient offline stage for PMOR.

\item As will be explained later, the choice of the number of basis vectors $n$ in the projection matrix $\bV$ can be determined through an adaptively evolving criterion calculated from the error estimator. This ensures just the adequate amount of basis vectors to guarantee the desired approximation quality with the ROM.
\end{itemize}

\subsection{An adaptive POD-Greedy algorithm to compute the basis $\bV$}
The POD-Greedy algorithm (\textsf{PODg}) was introduced in \cite{morHaaO08} as an extension of the greedy algorithm~\cite{morQuaMN16,morHesRS16} to time-dependent systems. Over the years, it has evolved to be the work-horse of the reduced basis method applied to time-dependent systems in a range of applications. In our previous works~\cite{morCheFB19a,morCheFB22}, we proposed extensions of the POD-Greedy algorithm, called adaptive POD-Greedy-DEIM algorithm (\textsf{aPODg+EI}), with several enhancements over the state-of-the-art. In this work, we make use of \textsf{aPODg+EI} to compute suitable ROMs for the cardiac electrophysiology equations~\cref{eq:fullydisc_monodomain_coupled}. We briefly review the main features of the algorithm and refer to~\cite{morCheFB19a,morCheFB22} for the details. 

The adaptive POD-Greedy algorithm constructs the projection basis $\bV$ through an iterative procedure, driven by an a posteriori output error estimator. The pseudo-code is sketched in~\Cref{alg:aPODg+EI}. The inputs to the algorithm are the ROM tolerance $\texttt{tol}$, a fine discretization of the parameter space $\mathcal{P}$ in the form of a training set $\Xi = \{ \p_{1}, \p_{2}, \ldots, \p_{n_{s}}\}$ and the FOM system matrices~(see \cref{eq:fullydisc_monodomain_coupled}). The outputs of the algorithm are the projection matrix $\bV$ and the hyperreduction quantities ($\bU, \bP$).

At any given iteration, the FOM~\cref{eq:fullydisc_monodomain_coupled} is solved at the current greedy parameter (denoted $\p^{*}$) and the resulting snapshot matrix $\bX_{s}(\p^*) \in \R^{N \times N_{t}}$ is used to enrich the projection basis $\bV$ with $n$ new basis vectors. The solution snapshots are used to compute the nonlinear snapshot matrix $\bF_{s}(\p^{*}) \in \R^{N \times N_{t}}$, which is used by the chosen hyperreduction algorithm (DEIM or EIM) to update $(\bU, \bP)$ with $n_{\text{EI}}$ basis vectors and selection indices. The greedy parameter for the subsequent iteration is chosen using the error estimator $\Delta(\p)$ based on the optimization problem
\[
	\p^* := \arg \max \limits_{\p \in \Xi} \Delta(\p).
\] 
The error estimator $\Deltamu$ can be evaluated once the ROM solutions~\cref{eq:fullydisc_monodomain_coupled_romEI} are computed. The greedy algorithm is said to have converged, when the maximum estimated error 
\[\epsilon := \Delta(\p^*) = \DeltamustarRB + \DeltamustarEI\]
is below the desired tolerance, i.e., $\epsilon < \texttt{tol}$. Here, $\DeltamustarRB$ is the maximum estimated error contributed by the reduced basis approximation and $\DeltamustarEI$ is the contribution to the estimated error due to hyperreduction. The number of basis vectors $n, \nEI$ is based on the update rule~\cite{morCheFB19a} determined from the estimated error
\begin{subequations}
\label{eq:basis_update_rules}
\begin{align}
\label{eq:basis_update_rules_RB}
n &:= \cRB \left\lfloor\log_{10}\left(\frac{\DeltamustarRB}{\texttt{tol}}\right)\right\rfloor,\\[0.8em]
\label{eq:basis_update_rules_EI}
\nEI &:= \nEI + \cEI \left\lfloor\log_{10}\left( \frac{\DeltamustarEI}{\texttt{tol}}\right)\right\rfloor.
\end{align}
\end{subequations}
\begin{remark}
\label{rem:basisV_computation}
Note that when the POD-Greedy algorithm is applied to coupled problems, e.g.,~\cref{eq:fullydisc_monodomain_coupled}, the solution snapshots corresponding to each variable, viz., $\bX_{s}^{\Phi}, \bX_{s}^{r} \in \R^{N \times N_{t}}$ are collected and two separate SVDs are performed to obtain the bases $\bV^{\Phi}, \bV^{r}$. Using these, the projection basis $\bV$ is formed as in~\cref{eq:blockdiag_V}.
\end{remark}
\begin{remark}
The rough intuition behind the update rule~\cref{eq:basis_update_rules} is that we add one basis vector per order of magnitude difference from the desired tolerance; e.g., suppose $\DeltamustarRB = 10^{-1}$ and $\texttt{tol} = 10^{-6}$, $n = 5$ implying that, potentially, if $5$ new basis vectors are added in the next iteration to the basis $\bV$, the estimated error $\Delta(\p^{*})$ is expected to decrease below the tolerance. But, based on the application, it is not necessary that one basis vector leads to one order of magnitude reduction of the estimated error. In the original version of the above update rule in~\cite{morCheFB19a}, $\cRB, \cEI$ were set to $1$. In this work, we have generalized this, such that $\cRB, \cEI \in \mathbb{Z}_{+}$ can be any positive integer. This is especially desirable for problems whose singular value decays slowly, e.g., convection-dominated problems and problems with travelling shocks such as in cardiac electrophysiology. $\cRB$ and $\cEI$ are hyperparameters; they can either be fixed or determined heuristically based on the singular value decay of the current greedy snapshot matrix $\bX_{s}(\p^*)$. 
\end{remark}
\begin{algorithm}[t!]
\caption{Adaptive POD-Greedy-(D)EIM (\textsf{aPODg+EI}) algorithm}
\begin{algorithmic}[1]
\label{alg:aPODg+EI}
\REQUIRE Training set $\Xi$, tolerance (\texttt{tol}), discretized system matrices $(\bbE, \bbA, \bM_{f}, \bB, \bff, \bx_{0})$

\ENSURE $\mathbf{V}$, DEIM matrices $(\bU, \bP)$
\vspace{0.5em}

{Initialize: $\bV = [\,]$, $n = n_{0}$, $\nEI = \nEIo$, $\bU = [\,], \bP = [\,]$, $\bF = [\,]$, randomly selected initial greedy parameter $\p^* \in \Xi$, $\epsilon = 1 + \texttt{tol}$.}
\vspace{0.5em}

\STATE Compute dual system projection matrix $\bVdu$ needed for the output error estimation

\WHILE{$\epsilon$ $>$ \texttt{tol}}
\vspace{0.5em}
	
\STATE Obtain FOM~\cref{eq:fullydisc_monodomain_coupled} snapshots $\bX_{s}(\p^*)$ at greedy parameter; \\compute nonlinear snapshots $\bF_{s}(\p^*)$
\vspace{0.5em}

\STATE Determine $\bU_{\bX}$ through the SVD of $\bar{\bX} := \bX_{s}(\p^*) - \bV \bV^{T} \bX_{s}(\p^*)$, where $\bU_{\bX}$ is the matrix of left singular vectors of $\bX_{s}(\p^*)$
\vspace{0.5em}

\STATE Update $\bV$ as $\bV := \texttt{orth}\left(\left[\bV,\,\bU_{\bX}(:,\, 1:n)\right]\right)$\\ with $\texttt{orth}\big(\cdot\big)$ denoting an orthogonalization process which can be implemented using the modified Gram-Schmidt process, or QR factorization
\vspace{0.5em}

\STATE Form nonlinear snapshot matrix $\bF := \left[\bF,\,\bF_{s}(\p^*)\right]$ and apply DEIM algorithm to get updated hyperreduction quantities $(\bU, \bP)$ with $\nEI$ basis vectors and interpolation indices
\vspace{0.5em}

\STATE Obtain reduced system matrices through Galerkin projection; solve ROM~\cref{eq:fullydisc_monodomain_coupled_romEI} to compute error estimator $\Deltamu$ $\forall \p \in \Xi$
\vspace{0.5em}

\STATE Select greedy parameter $\p^* := \arg \max \limits_{\p \in \Xi} \Delta(\p)$; set $\epsilon = \Delta(\p^*)$
\vspace{0.5em}

\STATE Compute the basis updates $n, \nEI$ for the next iteration based on the update rules~\cref{eq:basis_update_rules}
\vspace{0.5em}
\ENDWHILE
\end{algorithmic}
\end{algorithm}
\subsection{A posteriori output error estimation}
\label{subsec:apostouterr}
Accurate estimation of the error incurred by the ROM is critical for the success of the \textsf{aPODg+EI} algorithm. For the cardiac electrophysiology model, the quality of approximation of the output quantities (ECG or flux) is of particular interest, as these are the ones of medical consequence. A posteriori error estimation for output quantities has received considerable attention in the reduced basis community~\cite{morGre05,morHaaO11,morZhaFLetal15,morCheFB19a,morFenCB23}. We use the residual-based primal-dual a posteriori output error estimator proposed in~\cite{morCheFB19a}. So far, error estimation for cardiac electrophysiology has not been widely discussed in the literature.

The output error estimator $\Delta^{k}(\p)$ at a given parameter $\p$ and at a given time instance $t^{k}$ has the form
\begin{align}
\label{eq:outputerrestm}
\|\by^{k}(\p) - \bar{\by}^{k}(\p)\| \lessapprox \bigg(\bar{\rho}\, \beta \,\| \brdu \| + \lvert 1 - \bar\rho \rvert\, \| \btxdu \|\bigg)\, \| \br^{k} \| =: \Delta^{k}(\p).
\end{align}
Here, the residual resulting from the ROM~\cref{eq:fullydisc_monodomain_coupled_rom} with respect to the FOM~\cref{eq:fullydisc_monodomain_coupled} at the time instance $t^{k}$ is
\begin{align}
\label{eq:rpr}
\br^{k} := \bbA \btx^{k-1}(\p) + \delta t \left( \bM_{f} \bff(\btx^{k-1}, \p) + \bB i_{s}^{k} \right) - \bbE \btx^{k}(\p).
\end{align}
In practice, when hyperreduction is enforced, the ROM~\cref{eq:fullydisc_monodomain_coupled_romEI} is solved. The residual can be additively decomposed into two parts: one relating to the error due to the reduced basis approximation and the other relating to the error resulting from hyperreduction
\begin{align*}
\br^{k} &= \bbA \btx^{k-1}(\p) + \delta t \left( \bM_{f} \bff(\btx^{k-1}, \p) + \bB i_{s}^{k} \right) - \bbE \btx^{k}(\p),\\

&= \bbA \btx^{k-1}(\p) + \delta t \left( \bM_{f} \bff(\btx^{k-1}, \p) + \bM_{f} \bff_{\text{EI}}(\btx^{k-1}, \p) - \bM_{f} \bff_{\text{EI}}(\btx^{k-1}, \p) + \bB i_{s}^{k} \right) - \bbE \btx^{k}(\p),\\

&= \underbrace{\bbA \btx^{k-1}(\p) + \delta t \left( \bM_{f} \bff_{\text{EI}}(\btx^{k-1}, \p) + \bB i_{s}^{k} \right) - \bbE \btx^{k}(\p)}_{\br^{k}_{\text{RB}}} + \underbrace{\delta t \bM_{f} \left( \bff(\btx^{k-1}, \p) - \bff_{\text{EI}}(\btx^{k-1}, \p)\right)}_{\br^{k}_{\text{EI}}}.

\end{align*}

The constant $\bar{\rho}$ is estimated at every greedy iteration based on the available snapshots in $\bX_{s}(\p^*)$. The details of its estimation may be found in~\cite{morCheFB19a}. The quantity $\beta$ corresponds to the inf-sup-constant and in the case of the matrix spectral norm, it corresponds to the inverse of the smallest singular value of the matrix $\bbE$, i.e., $\beta := \frac{1}{\sigma_{\min}(\bbE)}$.

A dual system is involved in obtaining the error estimator and it is given by
\begin{align}
\label{eq:dualsystem}
\bE_{\text{du}} \bxd = \bC_{\text{du}},
\end{align}
where $\bE_{\text{du}} := \bbE^{T} \in \R^{2N \times 2N}, \bC_{\text{du}} := -\bC^{T} \in \R^{2N}$. The ROM corresponding to the dual system is
\begin{align}
\label{eq:dualsystem_rom}
\bhE_{\text{du}} \bhxd = \bhC_{\text{du}},
\end{align}
where $\bhE_{\text{du}} := \bVdu^{T} \bE_{\text{du}} \bVdu \in \R^{2N \times n_{\text{du}}}, \bhC_{\text{du}} := \bVdu^{T} \bC_{\text{du}} \in \R^{n_{\text{du}}}$. The dual projection matrix is $\bVdu \in \R^{2N \times n_{\text{du}}}$. As the dual system is not parameter-dependent, the dual  basis $\bVdu$ is obtained by applying a Krylov subspace method to the dual FOM~\cref{eq:dualsystem} as done in~\cite{morCheFB19a}.

Based on the dual FOM and ROM~(\cref{eq:dualsystem,eq:dualsystem_rom}), the dual residual has the form
\begin{align}
\label{eq:rdu}
\brdu := \bC_{\text{du}} - \bE_{\text{du}} \btxdu.
\end{align}
\paragraph{Mean estimated error}
The mean value of the estimated error over time, at a given parameter sample $\p$ is
\[
	\Deltamu := \frac{1}{N_{t}} \sum\limits_{k=0}^{K} \Delta^{k}(\p).
\]

\subsection{Adaptive choice of the training set}
\label{subsec:adapt_trngset}
So far, we have seen how to derive a projection-based ROM for the discretized monodomain equations for cardiac electrophysiology~\cref{eq:disc_monodomain_coupled} and discussed an adaptive algorithm \textsf{aPODg+EI}~(\Cref{alg:aPODg+EI}) to obtain the projection basis for the ROM~$\bV$. The \textsf{aPODg+EI} algorithm calls for an efficient error estimator to drive the greedy parameter sampling; such an error estimator for the output quantity is discussed in the previous section. However, the choice of the training set to be used in~\Cref{alg:aPODg+EI} is still unclear. In this section, we show how adaptivity can be used to systematically update the training set by adding (or removing) parameter samples to (from) it. Our approach for the adaptive sampling of the training set is based on the method proposed in~\cite{morCheFB22}. We briefly review this method next and refer the reader to~\cite{morCheFB22} for the finer aspects and implementation details. The adaptive greedy algorithm implementing the adaptive training set sampling is sketched in \Cref{alg:aPODg+EI+adaptTS}. Note that it is a generalization of~\Cref{alg:aPODg+EI}.
\begin{algorithm}[t!]
\caption{Adaptive POD-Greedy-(D)EIM with adaptive training set sampling (\textsf{aPODg+EI+adaptTS}) algorithm}
\begin{algorithmic}[1]
\label{alg:aPODg+EI+adaptTS}
\REQUIRE Coarse training set $\Xi_{c}$, fine training set $\Xi_{f}$, tolerance (\texttt{tol}), discretized system matrices $(\bbE, \bbA, \bM_{f}, \bB, \bff, \bx_{0})$

\ENSURE $\mathbf{V}$, DEIM matrices $(\bU, \bP)$
\vspace{0.5em}

{Initialize: $\bV = [\,]$, $n = n_{0}$, $\nEI = \nEIo$, $\bU = [\,], \bP = [\,]$, $\bF = [\,]$, randomly selected initial greedy parameter $\p^* \in \Xi_{c}$, $\epsilon = 1 + \texttt{tol}$.}
\vspace{0.5em}

\STATE Compute dual system projection matrix $\bVdu$ needed for the output error estimation

\WHILE{$\epsilon$ $>$ \texttt{tol}}
\vspace{0.5em}
	
\STATE Obtain FOM~\cref{eq:fullydisc_monodomain_coupled} snapshots $\bX_{s}(\p^*)$ at greedy parameter; \\compute nonlinear snapshots $\bF_{s}(\p^*)$
\vspace{0.5em}

\STATE Determine $\bU_{\bX}$ through the SVD of $\bar{\bX} := \bX_{s}(\p^*) - \bV \bV^{T} \bX_{s}(\p^*)$, where $\bU_{\bX}$ is the matrix of left singular vectors of $\bX_{s}(\p^*)$
\vspace{0.5em}

\STATE Update $\bV$ as $\bV := \texttt{orth}\left(\left[\bV,\,\bU_{\bX}(:,\, 1:n)\right]\right)$\\ with $\texttt{orth}\big(\cdot\big)$ denoting an orthogonalization process which can be implemented using the modified Gram-Schmidt process, or QR factorization
\vspace{0.5em}

\STATE Form nonlinear snapshot matrix $\bF := \left[\bF,\,\bF_{s}(\p^*)\right]$ and apply DEIM algorithm to get updated hyperreduction quantities $(\bU, \bP)$ with $\nEI$ basis vectors and interpolation indices
\vspace{0.5em}

\STATE Obtain reduced system matrices through Galerkin projection; solve ROM~\cref{eq:fullydisc_monodomain_coupled_romEI} to compute error estimator $\Deltamu$ $\forall \p \in \Xi_{c}$
\vspace{0.5em}

\STATE Obtain the radial basis interpolant and evaluate the (interpolated) error $\widetilde{\Delta}(\p)$ $\forall \p \in \Xi_{f}$
\vspace{0.5em}

\STATE Remove parameter samples $\mathring{\p} \in \Xi_{c}$ for which $\Deltamu < \texttt{tol}$ $\forall \p \in \Xi_{c}$
\vspace{0.5em}

\STATE Update $\Xi_{c}$ with $n_{\text{add}}$ new parameter samples taken from $\Xi_{f}$ with the highest evaluated error $\widetilde{\Delta}(\p)$
\vspace{0.5em}

\STATE Select greedy parameter $\p^* := \arg \max \limits_{\p \in \Xi} \Delta(\p)$; set $\epsilon = \Delta(\p^*)$
\vspace{0.5em}

\STATE Compute the basis updates $n, \nEI$ for the next iteration based on the update rules~\cref{eq:basis_update_rules}
\vspace{0.5em}
\ENDWHILE
\end{algorithmic}
\end{algorithm}

The proposed adaptive training set sampling method makes use of two training sets - a coarse training set denoted $\Xi_{c}$, consisting of $n_{c}$ parameter samples and a fine training set $\Xi_{f}$, having $n_{f} \gg n_{c}$ parameter samples. The error estimator $\Deltamu$ is evaluated only for the samples present in $\Xi_{c}$. We learn a radial basis interpolant of the mapping $\Delta: \p \xrightarrow{} \Deltamu$ $\forall \p \in \Xi_{c}$. Following this, the interpolant is queried to evaluate the error at the parameter samples in $\Xi_{f}$, i.e., $\widetilde{\Delta}(\p) = \chi(\p)$ $\forall \p \in \Xi_{f}$ (Step 8 in~\Cref{alg:aPODg+EI+adaptTS}). Here, $\chi(\cdot)$ is the radial basis interpolant function. Computational efficiency comes from the fact that, at every greedy iteration, the ROM needs to be evaluated only at the parameter samples in the coarse training set. As illustrated in~\cite{morChe23}, constructing and evaluating the radial basis interpolant over the fine training set has comparatively less cost than solving the ROM at all parameter samples in the fine training set.

At the end of each greedy iteration, the coarse training set $\Xi_{c}$ is updated with new parameter samples taken from $\Xi_{f}$. We add $n_{\text{add}}$ new parameter samples that have the largest error evaluated by the interpolant $\forall \p \in \Xi_{f}$. In addition to this, we also examine the coarse training set to identify and remove samples $\mathring{\p} \in \Xi_{c}$ for which $\Delta(\mathring{\p}) < \texttt{tol}$. Doing this two-pronged approach of adding and removing samples from $\Xi_{c}$ ensures that the training set remains as compact as possible.
\section{Numerical experiments}%
\begin{figure}[t!]
\centering
\includegraphics[scale=1]{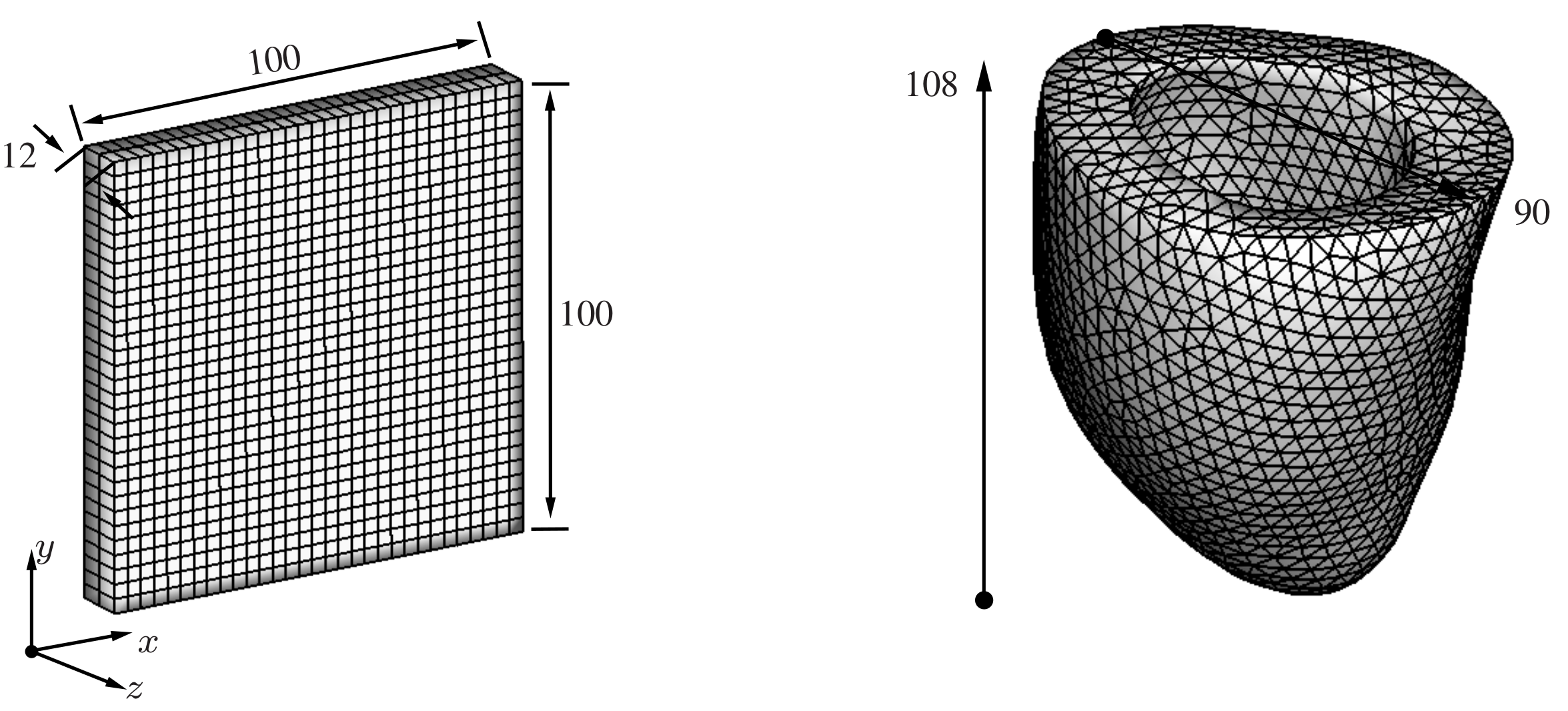}
\caption{Discretization of \textbf{Left:} the myocardial tissue slice and \textbf{Right:} the left ventricle model of a healthy person. All lengths are in millimetres.}
\label{fig:3dblock}
\end{figure}

\label{sec:numerics}
We apply two adaptive algorithms~(\Cref{alg:aPODg+EI} and \Cref{alg:aPODg+EI+adaptTS}) to two different benchmark examples which are
\begin{enumerate}
\item 3-D cardiac tissue block, see \Cref{fig:3dblock} (left), 
\item Personalized left ventricle (LV) model, see \Cref{fig:3dblock} (right).
\end{enumerate}
In both cases, the quantity of interest is the summation of electrical flux over the domain computed as $\displaystyle q = \int_{\mathcal{B}}\mathbf{q} \cdot\mathbf{n}~d\mathbf{z}$ with $\mathbf{n}$ denoting the projection direction of the flux vector $\mathbf{q} \in \R^{3}$. 
The projection direction $\mathbf{n}$ coincides with the $x$-axis in the tissue block examples while for the LV model, the longitudinal axis is considered as the projection direction. Note that the plot of $q$ versus time in the LV model corresponds to ECG. To verify the quality of the ROMs resulting from the application of our proposed adaptive algorithms, we use the following metrics

\begin{itemize}
\item The \emph{scaled} maximal estimated error at each greedy iteration is denoted by $\epsilon_{\text{max}}$ which is defined as
\[
	\epsilon_{\text{max}} := \max \limits_{\p \in \Xi_{\text{train}}} \frac{\Delta(\p)}{\texttt{scaling}},
\]
where $\texttt{scaling}$ is a scaling factor pre-defined by the user. Different options for the $\texttt{scaling}$ are possible. One can use a maximal scaling where at each iteration (of~\Cref{alg:aPODg+EI} or \Cref{alg:aPODg+EI+adaptTS}), the scaling factor is the maximum of the output evaluated at the current greedy parameter, i.e., $\texttt{max}\left(\bY(\p^{*})\right)$. Another potential option would be to use $\texttt{scaling} = \| \bY(\p^{*}) \|_{2}$. A further approach would be to use the average of $\| \bY(\p_{i}^{*}) \|_{2}$ where $\p_{i}^{*}$ are the already selected parameters during the greedy iterations. We use the maximum norm scaling, viz., $\texttt{max}\left(\bY(\p^{*})\right)$ in the presented numerical results.

\item The relative true error at a given parameter $\p$ as denoted by $\epsilon_{\text{rel}}$ is
\[
	\epsilon_{\text{rel}}(\p) := \frac{\| \bY(\p) - \widehat{\bY}(\p) \|_{2}}{\| \bY(\p) \|_{2}},
\]
where $\bY := \big[\by^{0}, \by^{1}, \ldots, \by^{K}\big]^{T} \in \R^{N_{t}}$ and $\widehat{\bY} := \big[\bhy^{0}, \bhy^{1}, \ldots, \bhy^{K}\big]^{T} \in \R^{N_{t}}$.
\end{itemize}

We refer by \textsf{Test A} the application of~\Cref{alg:aPODg+EI} to an example and by \textsf{Test B} the application of~\Cref{alg:aPODg+EI+adaptTS} to an example. 

\subsection{3-D cardiac tissue block}
In this example, we consider the domain of interest to be a 3-D block possessing the material properties of the cardiac tissue. The domain is discretized by 31x31x2 brick elements over 3072 nodes. We are interested in the variation of ECG-like diagrams as a function of the parameter $\gamma$ (appearing in the Aliev-Panfilov model~(see \cref{eq:AP_Fr})) and time, i.e., $\mathbf{q} := \bD \cdot \nabla \Phi(t, \gamma)$. To this end, the range of parameter variations is $\p := \gamma \in \mathcal{P} := [0.0005,\, 0.01]$. In the following, two scenarios are considered: a regular planar wave propagation and scroll wave propagation. For the case of scroll wave generation, a second parameter will be considered. 

\subsubsection{3-D cardiac tissue block with planar wave}
For the first example, where a regular wave propagation is considered, we model the conduction of the cardiac action potential initiated through an initial condition $\bx_{0} = \big[\bx_{0}^{\Phi} \, \bx_{0}^{r}\big]^{T} \in \R^{2\cdot3072}$ applied to the left boundary of the 3-D block. In order to initiate the planar wave propagation, the initial conditions of the nodes on the left edge are set to $\Phi_0= -10$ mV and for the remaining nodes $\Phi_0=-80$ mV. The initial condition for the recovery variable is set to zero at all nodes, i.e., $r_0 = \mathbf{0}$. The time step is set to $\delta t = 2$ milliseconds.

\begin{figure}
\centering
\includegraphics[scale=0.8]{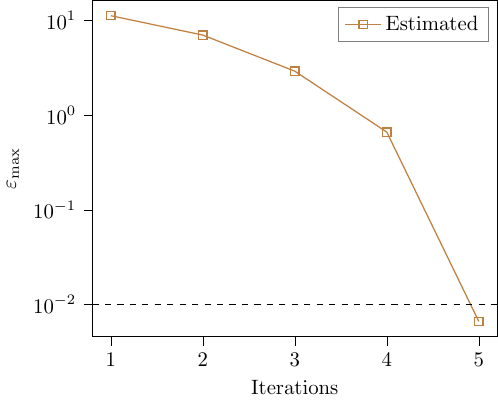}
\hfill
\includegraphics[scale=0.8]{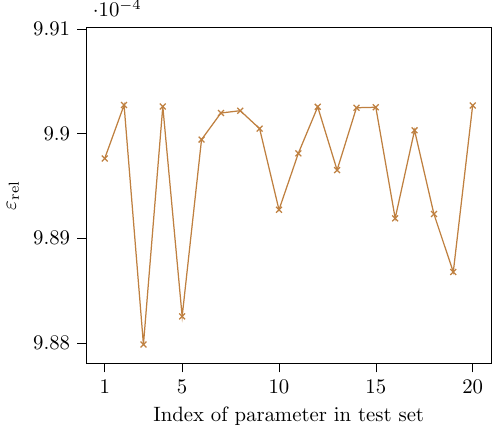}
\caption{\Cref{alg:aPODg+EI} applied to 3-D block of cardiac tissue: \textbf{Left:} Convergence of the greedy algorithm; \textbf{Right:} Performance of ROM on test set.}
\label{fig:alg1_3dblock}
\end{figure}

\begin{figure}
\centering
\includegraphics[scale=0.8]{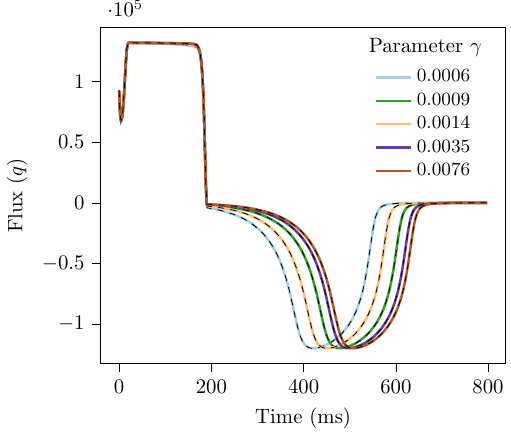}
\caption{\Cref{alg:aPODg+EI} applied to 3-D block of cardiac tissue: Comparison of flux obtained from the FOM (solid line) and the ROM (dashed line) at five different values of the parameter $\gamma$.}
\label{fig:alg1_3dblock_flux}
\end{figure}

\begin{figure}
\centering
\includegraphics[scale=0.8]{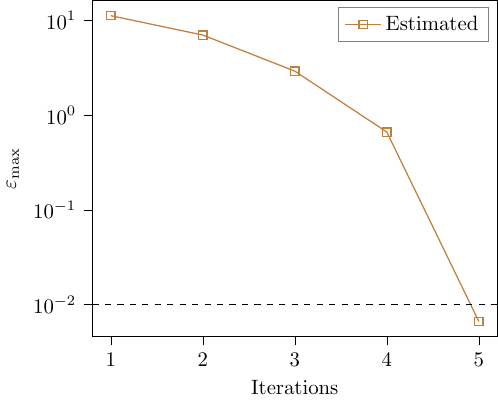}
\hfill
\includegraphics[scale=0.8]{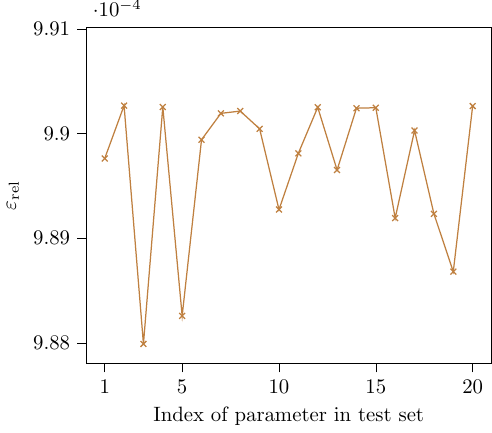}
\caption{\Cref{alg:aPODg+EI+adaptTS} applied to 3-D block of cardiac tissue: \textbf{Left:} Convergence of the greedy algorithm; \textbf{Right:} Performance of ROM on test set.}
\label{fig:alg2_3dblock}
\end{figure}

\begin{figure}
\centering
\includegraphics[scale=0.8]{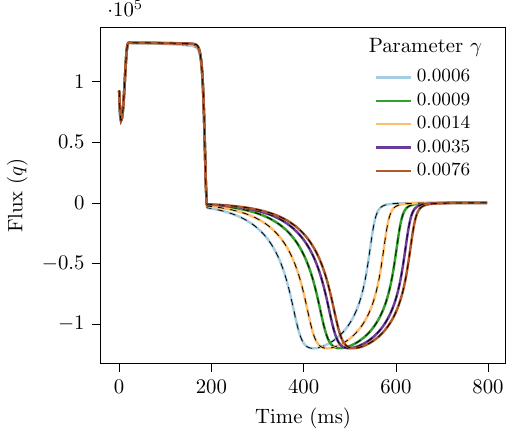}
\caption{\Cref{alg:aPODg+EI+adaptTS} applied to 3-D block of cardiac tissue: Comparison of flux obtained from the FOM (solid line) and the ROM (dashed line) at five different values of the parameter $\gamma$.}
\label{fig:alg2_3dblock_flux}
\end{figure}

\begin{figure}
\centering
\includegraphics[scale=0.8]{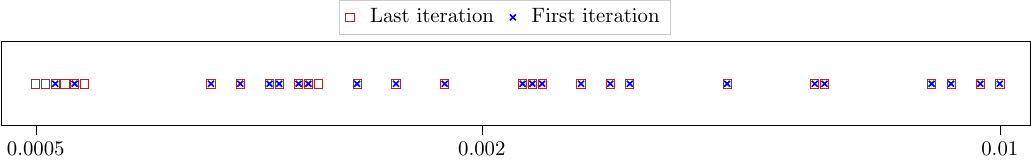}
\caption{\Cref{alg:aPODg+EI+adaptTS} applied to 3-D block of cardiac tissue: Comparison of the coarse training set in the first and last iteration.}
\label{fig:alg2_3dblock_trngset_first_last}
\end{figure}

\paragraph{\textsf{Test A} -- Application of \Cref{alg:aPODg+EI}:} We apply \Cref{alg:aPODg+EI} to the discretized monodomain equation in the 3-D block. To obtain a parameter set $\Xi$, we discretize the parameter domain $\mathcal{P}$ to collect $100$ parameter samples. We divide this set randomly in the ratio $80\,:\,20$ to obtain, respectively, a training set $\Xi_{\text{train}}$ and a test set $\Xi_{\text{test}}$. The tolerance $\texttt{tol} = 10^{-2}$ is used. We set $\cRB = \cEI = 1$. To determine the initial values of the reduced basis dimension, i.e., $n_{0}$, we make use of the singular value decay of the snapshots matrices $\bX_{s}^{\Phi}, \bX_{s}^{r}$ in the first iteration (see~\Cref{rem:basisV_computation}). The tolerances for the singular value decay in both cases are set to $0.5$. This results in $n_{0}^{\Phi} = 1, n_{0}^{r} = 2$. Further, we let $\nEIo^{\phi} = \nEIo^{r} = 8$. The greedy algorithm converges to the desired tolerance in $5$ iterations and requires $71$ seconds. The dimension of the resulting ROM is $n = 57$ with $n^{\Phi} = 28, n^{r} = 29$. The DEIM basis has dimension $\nEI = 64$ with $\nEI^{\phi} = \nEI^{r} = 32$. In terms of speedup achieved, while a single FOM solution requires $0.484$ seconds, a single ROM solution needs $0.016$ seconds, a $30$-fold acceleration. In~\Cref{fig:alg1_3dblock}, in the left figure, we plot the convergence of the maximum estimated error. On the right, the performance of the ROM obtained from~\Cref{alg:aPODg+EI} is plotted. We evaluate the ROM for the $20$ parameter samples in $\Xi_{\text{test}}$. It is evident that the relative true error $\epsilon_{\text{rel}}$ is below the desired tolerance for every sample. The flux waveform obtained from the FOM and the ROM for five different values of $\gamma$ are shown in~\Cref{fig:alg1_3dblock_flux}. Both waveforms are visually indistinguishable, showing their excellent agreement. We have thus reduced the dimension of the coupled system~\cref{eq:fullydisc_monodomain_coupled} from $N=2\cdot3072$ to $n=59$, while ensuring that the ROM is accurate up to a desired tolerance. In fact, the ROM displays excellent performance over unseen parameter samples, showing its ability to generalize well.
\paragraph{\textsf{Test B} -- Application of \Cref{alg:aPODg+EI+adaptTS}:} Next, we apply~\Cref{alg:aPODg+EI+adaptTS}. The training set and test set, $\Xi_{\text{train}}, \Xi_{\text{test}}$ are obtained in the same fashion as above. The coarse training set $\Xi_{c}$ and the fine training set $\Xi_{f}$ are obtained, respectively, by dividing the training set $\Xi_{\text{train}}$ in a $30\,:\,70$ ratio, resulting in the coarse training set with $n_{c} = 24$ samples and the fine training set with $n_{f} = 56$ samples. The choice of $n_{0}$ and $\nEIo$ are the same as in the previous case. We set $n_{\text{add}} = 1$. The greedy algorithm converges in $34.5$ seconds to the set tolerance, taking $5$ iterations. Note that the time taken in this case is roughly half that taken in case of \textsf{Test A}. The dimension of the projection basis $\bV$ is $n = 57$ ($n^{\Phi} = 28, n^{r} = 29$). Furthermore, the dimension of $\bU$ is $\nEI = 64$, as for the previous case. Since the dimension of the ROM is same as for \textsf{Test A}, the same speedup of $30$x is achieved also for this case. \Cref{fig:alg2_3dblock} illustrates the results of applying \Cref{alg:aPODg+EI+adaptTS} to the 3-D cardiac block. In the left figure, we see the convergence of the estimated error while the right figure shows the relative true error $\epsilon_{\text{rel}}$ on the samples in $\Xi_{\text{test}}$. Once again, we see that the error is well below the desired tolerance, as in the previous case with the fixed training set. We also see from~\Cref{fig:alg2_3dblock_flux} an excellent agreement of the flux plotted as a function of time obtained using the FOM and the ROM at different values of $\gamma$. It is thus reliably demonstrated that iteratively building the training set leads to halving the time required to obtain the ROM; at the same time, the resulting ROM also meets the accuracy defined by the user. The evolution of the coarse training set $\Xi_{c}$ is shown in~\Cref{fig:alg2_3dblock_trngset_first_last}. In the first iteration it contains $24$ parameter samples (blue crosses). At the final iteration, we observe that there are now $29$ parameter samples (brown squares). Note that several new values of the parameter $\gamma$ are added in the leftmost region where $\gamma$ attains smaller magnitude. 
 

\subsubsection{3-D cardiac tissue block with scroll wave}
In this part, a scroll wave or, in other words, a reentrant wave is considered which is one of the fundamental benchmark problems in cardiac electrophysiology \cite{goktepe+wong+etal10}. The reentrant wave propagation leads to a chaotic electrical activity in the heart and is often named as arrhythmia, where the pumping function of the heart is diminished or even stops \cite{keating+sanguinetti01}.

One can simulate the generation of scroll waves by applying an appropriately timed stimulus to the monodomain equation. The initial condition, viz., $\bx_{0}$ is slightly different when compared to the previous example where no scroll wave generation is considered. We define the initial condition for the transmembrane potential $\bx_{0}^{\Phi}$ to be uniformly $-80$ mV at all nodes. The initial condition for the recovery variable $\bx_{0}^{r}$ is the same and set to $\bx_{0}^{r} = \mathbf{0}$. To trigger the propagation of the planar wave, we apply a stimulus for the first 5 time steps (i.e., 10 milliseconds duration). For the scroll wave generation, we consider the parameter $\p =(\gamma, t_{s})$ belonging to a two-dimensional parameter space $\mathcal{P} :=  [0.0015,\, 0.002]\, \times\, [480,\, 500]$. Here, $\gamma$ is the conductivity as before and the new parameter $t_{s}$ denotes the initiation time of the stimulus to trigger the scroll wave. The stimulus is applied for a span of $20$ milliseconds, or $10$ time steps. We note that this example is particularly challenging due to the variety of (stiff) dynamics exhibited by the system as a result of the scroll wave formation. The slow propagation of the action potential wavefront across the domain results in a rather slow decay of the singular values. Therefore, the ROM (and hyperreduction) basis needs many basis vectors to accurately capture the true dynamics.

\begin{figure}[t!]
\centering
\includegraphics[scale=0.8]{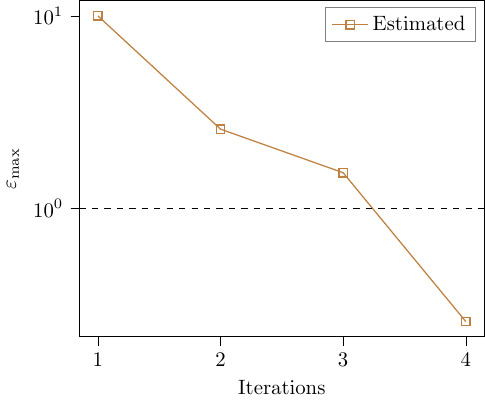}
\hfill
\includegraphics[scale=0.8]{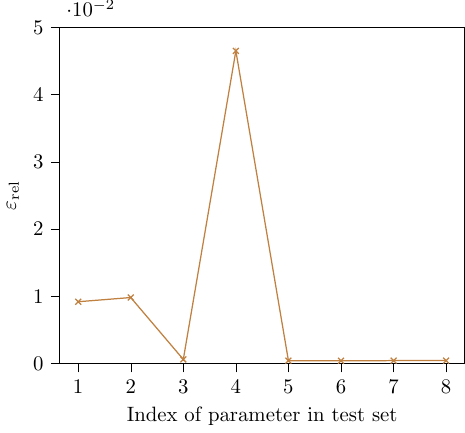}
\caption{\Cref{alg:aPODg+EI} applied to 3-D block of cardiac tissue with scroll wave generation: \textbf{Left:} Convergence of the greedy algorithm; \textbf{Right:} Performance of ROM on test set.}
\label{fig:alg1_3dblockscroll}
\end{figure}

\begin{figure}[t!]
\centering
\includegraphics[scale=0.8]{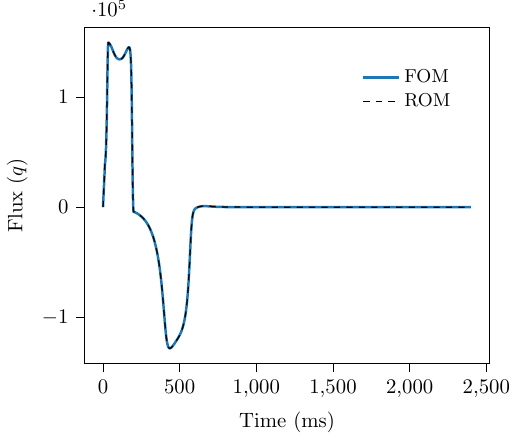}
\hfill
\includegraphics[scale=0.8]{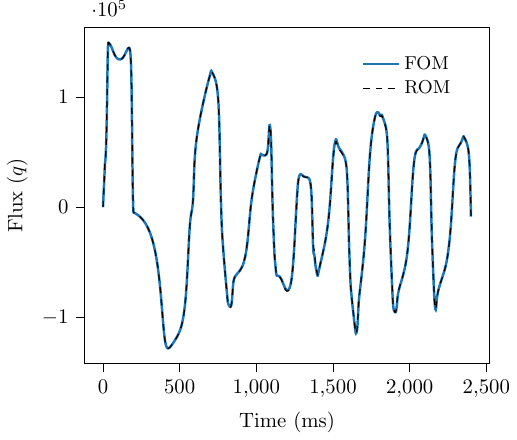}
\caption{\Cref{alg:aPODg+EI} applied to 3-D block of cardiac tissue with scroll wave generation: Comparison of flux obtained from the FOM (solid line) and the ROM (dashed line) at two different parameter samples $\p := (\gamma, t_{s})$ from the test set \textbf{Left:} Flux at $\p = (0.0017, 488)$; \textbf{Right:} Flux at $\p = (0.0020, 484)$.}
\label{fig:alg1_3dblockscroll_flux}
\end{figure}

\begin{figure}
\centering
\includegraphics[width=1\textwidth]{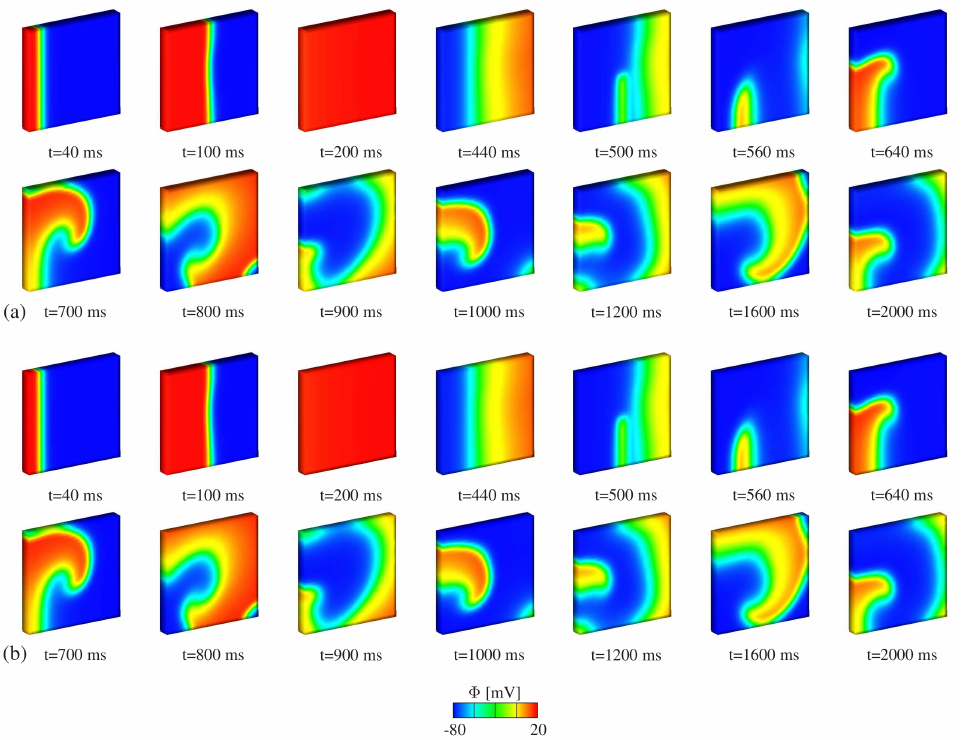}
\caption{Transmembrane potential $\Phi$ of the 3D block at different time instances: (a) FOM, (b) ROM from \textsf{Test B}.}
\label{fig:alg2_3dblock_scroll_Phi_FOM_ROM}
\end{figure}

\paragraph{\textsf{Test A} -- Application of \Cref{alg:aPODg+EI}:}
The tolerance of \Cref{alg:aPODg+EI} is set to be $\texttt{tol} = 1.0$ due to the fact that this example is challenging. The training set $\Xi$ is obtained by choosing $6 \times 6$ uniformly-spaced samples from the domain $\mathcal{P}$. This training set is divided in the ratio $80\, :\, 20$ to get $\Xi_{\text{train}}$ and $\Xi_{\text{test}}$. To select $n_{0}^{\Phi}, n_{0}^{r}$, the SVD tolerance is set to $10^{-5}$ and $10^{-3}$ in the first iteration, respectively, for the transmembrane potential snapshots $\bX_{s}^{\Phi}$ and the recovery voltage snapshots $\bX_{s}^{r}$. Owing to the fast changing dynamics of the nonlinearity, we determined the best setting for $\nEIo = (680 + 740) = 1420$. We further set $\cRB = 11, \cEI = 6$ as the factors to update the basis adaptively (see \cref{eq:basis_update_rules}). As seen from~\Cref{fig:alg1_3dblockscroll} (left), the greedy algorithm converges to the tolerance in $4$ iterations. The time taken is $1303$ seconds. The dimension of the basis $\bV$, i.e., $n = (572 + 264) = 836$ while the hyperreduction basis $\nEI = (778 + 830)  = 1608$. Evidently, this is a significantly high number of basis vectors and the reason for this is the slow decay of the singular values of the snapshot matrix. As a fallout of the large ROM dimension, no speedup is achieved for this test case. The FOM evaluation time is $1.298$ seconds, whereas the ROM evaluation time is slightly larger at $1.503$ seconds. Despite the larger ROM dimension, the performance on the test set $\Xi_{\text{test}}$ is satisfactory (see \Cref{fig:alg1_3dblockscroll} right figure). The relative error $\epsilon_{\text{rel}}$ is below the desired tolerance for all the test parameters. To illustrate the quality of the ROM approximation, we plot in~\Cref{fig:alg1_3dblockscroll_flux} the flux waveform for two different parameters in the test set. These are chosen such that the resulting flux exhibits different physical behaviour. For the first parameter $\p = (0.0017, 488)$, it can be seen that the scroll wave is not initiated. This is owing to the stimuli being applied late. The ROM is able to accurately capture this behaviour. In case of the second parameter $\p = (0.0020, 484)$, since the timing of the external stimuli is good, the scroll wave behaviour occurs. Once again, the ROM is also able to express this behaviour accurately.
\paragraph{\textsf{Test B} -- Application of \Cref{alg:aPODg+EI+adaptTS}:}
\begin{figure}[t!]
\centering
\includegraphics[scale=0.8]{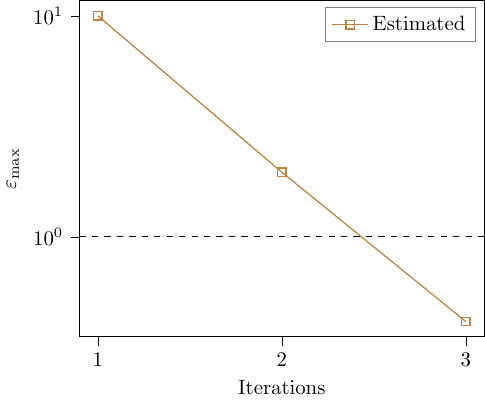}
\hfill
\includegraphics[scale=0.8]{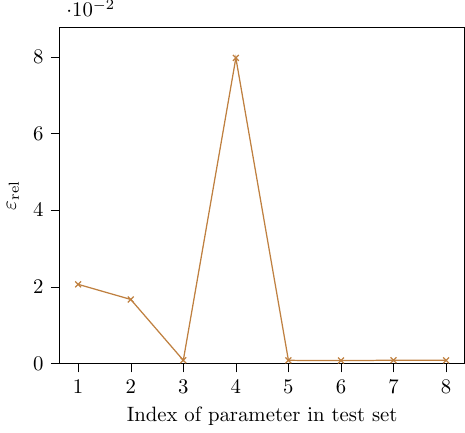}
\caption{\Cref{alg:aPODg+EI+adaptTS} applied to 3-D block of cardiac tissue with scroll wave generation: \textbf{Left:} Convergence of the greedy algorithm; \textbf{Right:} Performance of ROM on test set.}
\label{fig:alg2_3dblockscroll}
\end{figure}
\begin{figure}
\centering
\includegraphics[scale=0.8]{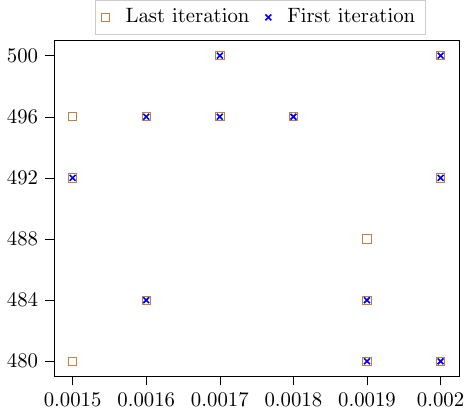}
\caption{\Cref{alg:aPODg+EI+adaptTS} applied to 3-D block of cardiac tissue with scroll wave generation: Comparison of the coarse training set in the first and last iteration.}
\label{fig:alg2_3dblock_scroll_trngset_first_last}
\end{figure}
\begin{figure}[t!]
\centering
\includegraphics[scale=0.8]{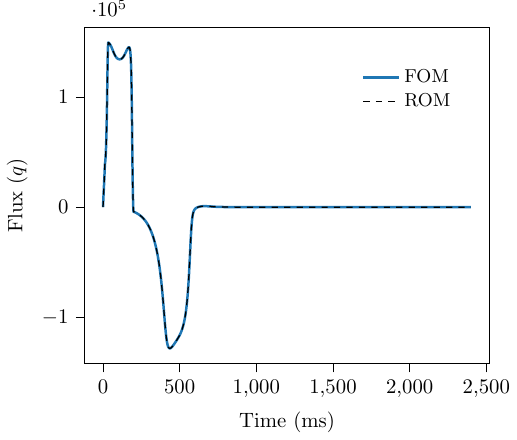}
\hfill
\includegraphics[scale=0.8]{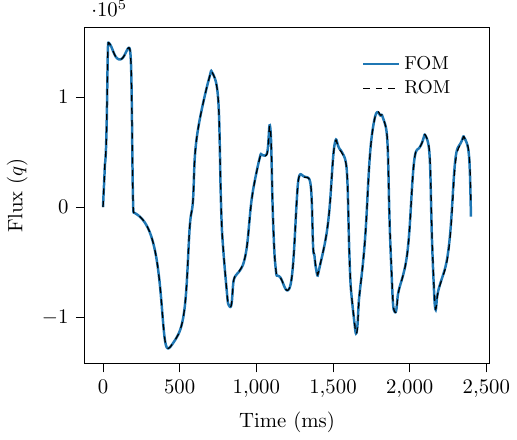}
\caption{\Cref{alg:aPODg+EI+adaptTS} applied to 3-D block of cardiac tissue with scroll wave generation: Comparison of flux obtained from the FOM (solid line) and the ROM (dashed line) at two different parameter samples $\p := (\gamma, t_{s})$ from the test set \textbf{Left:} Flux at $\p = (0.0017, 488)$; \textbf{Right:} Flux at $\p = (0.0020, 484)$.}
\label{fig:alg2_3dblockscroll_flux}
\end{figure}
Now we combine the adaptive greedy algorithm with the adaptive parameter sampling for the scroll wave example. The set $\Xi$ is same as before. We divide $\Xi$ in the ratio $40\, :\, 60$ to form the coarse training set and fine training set, respectively. The settings for $n_{0}$ and $\nEIo$ are retained as previously done. On the left figure of \Cref{fig:alg2_3dblockscroll}, the convergence of the greedy algorithm is plotted. We see that, in comparison to \textsf{Test A} we need only $3$ iterations. Correspondingly, \Cref{alg:aPODg+EI+adaptTS} takes only $492$ seconds to achieve convergence. The dimension of the ROM basis is smaller with $n = (517 + 209) = 726$ basis vectors in $\bV$; while the basis $\bU$ consists of $\nEI = (752 + 806) = 1558$ basis vectors. Observe that the ROM dimension in this test is smaller. As a result, we obtain a modest speedup. One FOM simulation takes $1.35$ seconds whereas the ROM needs $1.235$ seconds, a speedup of $1.1$x. The training set at the first and final iterations are shown in~\Cref{fig:alg2_3dblock_scroll_trngset_first_last}. The ROM obtained from this test also performs well on the test set $\Xi_{\text{test}}$. However, it is worth noting that the maximum error incurred on the test set is slightly higher for the current test (0.0797) while it is a little less for \textsf{Test A} (0.0465). This can be explained by the smaller ROM dimension in case of \textsf{Test B}. Nevertheless, as seen from~~\Cref{fig:alg2_3dblockscroll_flux}, the ROM obtained from \textsf{Test B} faithfully approximates the flux behaviour for two different parameter choices from the test set ($\p = (0.0017, 488)$ and $\p = $(0.0020, 484)) which exhibit different flux waveform patterns. We further plot the evolution of the transmembrane potential $\Phi$ at time instances $t_{k} \in \{40, 100, 200, 440, 500, 560, 640, 700, 800, 900, 1000, 1200, 1600, 2000\}$ milliseconds in~\Cref{fig:alg2_3dblock_scroll_Phi_FOM_ROM}. The $\Phi$ snapshots coming from the FOM and the ROM obtained with \textsf{Test B} show excellent agreement, both qualitatively and quantitatively.
\begin{remark}
We emphasize that the scroll wave propagation is a particularly challenging example. While both ROMs (from \textsf{Test A} and \textsf{Test B}) exhibit good performance over the test set, they do not offer significant speed up. This is owing to the larger ROM size. Linear projection-based MOR methods such as the ones used in this work are not very efficient for convection-dominated problems or problems with travelling waves~\cite{GreU19,morCagMS19}. Such problems typically require many basis vectors in the projection basis $\bV$ owing to the slow decay of the singular values. In case of cardiac electrophysiology, the nonlinearity and the coupled nature of the problem adds to the difficulty faced by linear MOR methods. Recent works have sought to address this using quadratic or even nonlinear manifolds. In addition, machine learning approaches have also been used. Incorporating these within our adaptive greedy framework to improve the speed up offered by ROMs will be a subject for future investigation.     
\end{remark}
\subsection{Left ventricle of the human heart}
Next, the introduced algorithms are tested on a personalized LV geometry that is generated from 4D echocardiography data of a healthy subject, see \cite{cansiz+sveric+etal18} for the procedure of model generation. The LV geometry is created at enddiastole and discretized by $19096$ four-node tetrahedral elements over $4129$ nodes. The total number of degrees of freedom in the large-scale FOM~\cref{eq:fullydisc_monodomain_coupled} is $N=2 \cdot 4129$. The parameter of interest is $\p := \gamma \in \mathcal{P} := [0.0002,\, 0.01]$. Moreover, the scaling factor $\beta_t$ is scaled down in terms of the activation time $t^{\text{act}}$ thereby early excited regions have longer action potential duration to the later activated regions as suggested in \cite{kotikanyadanam+goktepe+etal10}. The activation time $t^{\text{act}}$ is determined as the elapsed time from the atriventricular node stimulation until the transmembrane potential value of a particular point reaches -40 mV.

\begin{figure}
\centering
\includegraphics[width=1\textwidth]{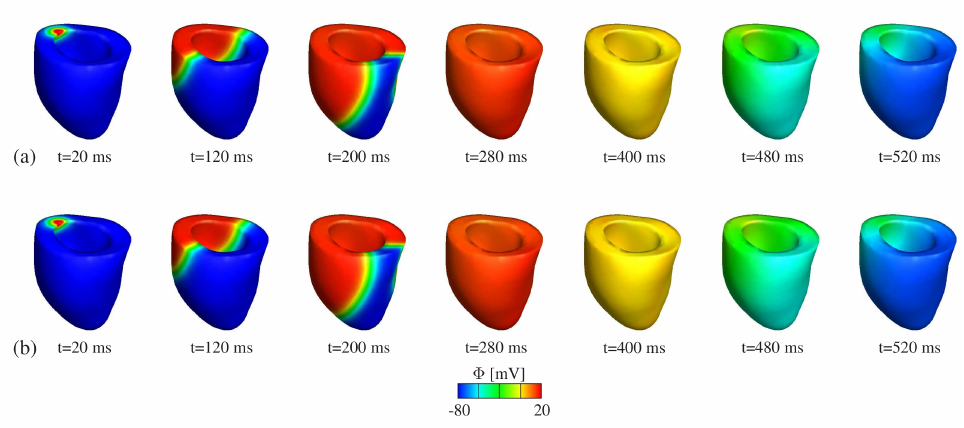}
\caption{Transmembrane potential $\Phi$ of the left ventricle at different time instances: (a) FOM, (b) ROM from \textsf{Test B}.}
\label{fig:alg2_ventricle_Phi_FOM_ROM}
\end{figure}

Similar to our previous works \cite{cansiz+dal+etal17,cansiz+sveric+etal18}, before the final results are obtained, three cycles are performed in order to achieve saturated values of the primary field variables over the subsequent cycles. Each cardiac cycle, which has a duration of 800 ms, is
initiated by applying a small amount of stimulus (I = 10 [-]) for 10 milliseconds to the upper part of the septum corresponding to the atrioventricular node (see the snapshots at time $t=20$ ms in~\Cref{fig:alg2_ventricle_Phi_FOM_ROM}).

\begin{figure}[t!]
\centering
\includegraphics[scale=0.8]{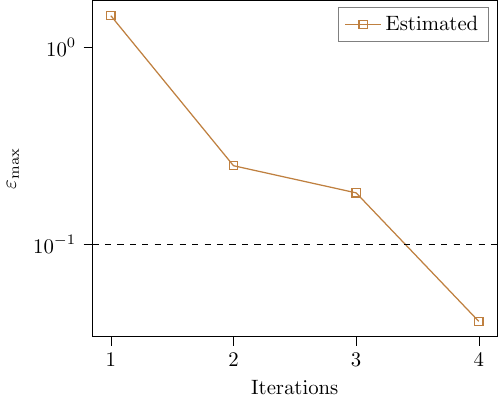}
\hfill
\includegraphics[scale=0.8]{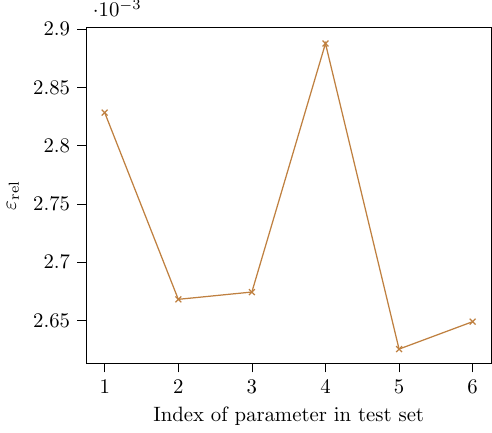}
\caption{\Cref{alg:aPODg+EI} applied to the left ventricle model: \textbf{Left:} Convergence of the greedy algorithm; \textbf{Right:} Performance of ROM on test set.}
\label{fig:alg1_ventricle}
\end{figure}

\begin{figure}[t!]
\centering
\includegraphics[scale=0.8]{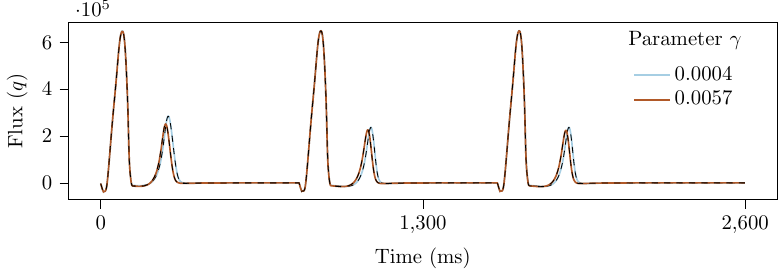}
\caption{\Cref{alg:aPODg+EI} applied to the left ventricle model: Comparison of flux obtained from the FOM (solid line) and the ROM (dashed line) at two different values of the parameter $\gamma$.}
\label{fig:alg1_ventricle_flux}
\end{figure}

\begin{figure}[t!]
\centering
\includegraphics[scale=0.8]{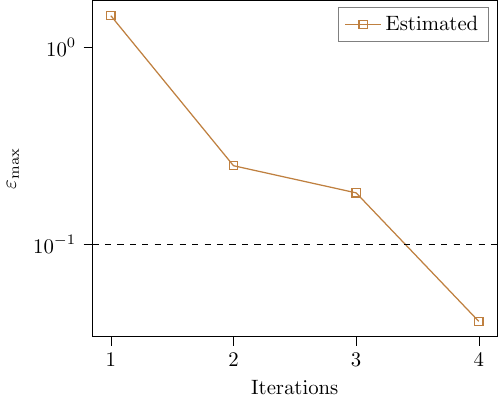}
\hfill
\includegraphics[scale=0.8]{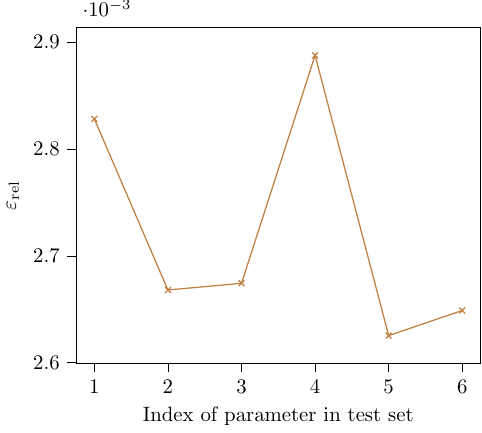}
\caption{\Cref{alg:aPODg+EI+adaptTS} applied to to the left ventricle: \textbf{Left:} Convergence of the greedy algorithm; \textbf{Right:} Performance of ROM on test set.}
\label{fig:alg2_ventricle}
\end{figure}

\begin{figure}[t!]
\centering
\includegraphics[scale=0.8]{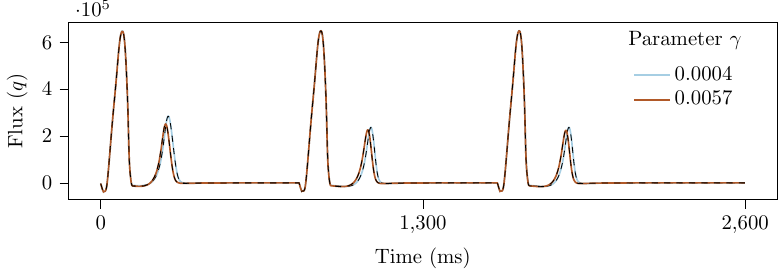}
\caption{\Cref{alg:aPODg+EI+adaptTS} applied to left ventricle model: Comparison of flux obtained from the FOM (solid line) and the ROM (dashed line) at two different values of the parameter $\gamma$.}
\label{fig:alg2_ventricle_flux}
\end{figure}

\paragraph{\textsf{Test A} -- Application of \Cref{alg:aPODg+EI}:}
For the LV geometry, we start by applying~\Cref{alg:aPODg+EI} to obtain a ROM. The training set $\Xi$ consists of $54$ samples of the parameter $\gamma$ obtained by selecting a pool of $60$ uniformly-spaced samples in the range $[0.0002,\, 0.01]$. It is further divided in the ratio $90\, :\, 10$ to form the training set $\Xi_{\text{train}}$ and the test set $\Xi_{\text{test}}$. In this example, we choose $n_{0}^{\Phi}$ and $n_{0}^{r}$ based on the SVD tolerances of $10^{-2}$; further, we set $\nEIo^{\phi} = 30$, $\nEIo^{r} = 10$. The convergence of the greedy algorithm is shown in the left figure in~\Cref{fig:alg1_ventricle}. The estimated error reaches the desired tolerance of $\texttt{tol} = 0.1$ in $4$ iterations and the time taken is $987$ seconds. The ROM dimension is $n = (40 + 36) = 76$ while the hyperreduction basis has size $\nEI = (113 + 74) = 187$. The accuracy of the obtained ROM is tested by predicting the flux values at the test parameter samples. It is clear from~\Cref{fig:alg1_ventricle} (right figure) that for all the samples, the relative error is less than $0.2 \%$ showing the good quality of approximation offered by the reduced model. In~\Cref{fig:alg1_ventricle_flux}, the flux resulting from the FOM and ROM simulations are plotted at two different values of $\gamma$ ($0.0004$, $0.0057$). The two waveforms display a good match. The time taken to solve the FOM is $6.646$ seconds. The ROM offers a speedup of around $49$x, taking only $0.135$ seconds to compute.
\paragraph{\textsf{Test B} -- Application of \Cref{alg:aPODg+EI+adaptTS}:}
Next, we show the benefits of the adaptive sampling of the training set for the LV model using~\Cref{alg:aPODg+EI+adaptTS}. The training set $\Xi_{\text{train}}$, the test set $\Xi_{\text{test}}$, the initialization for $n_{0}$ and $\nEIo$ follow the same values as those for \textsf{Test A}. The coarse training set $\Xi_{c}$ is obtained by choosing $30\%$ of the samples from $\Xi_{\text{train}}$ with $\Xi_{f}$ containing the remaining $70 \%$. We take $n_{\text{add}} = 1$. The tolerance is the same as before with $\texttt{tol} = 0.1$. The results are illustrated in~\Cref{fig:alg2_ventricle}. As seen in the left figure, the greedy algorithm converges in $4$ iterations. However, the time taken is only $507$ seconds, which is roughly half the time used in \textsf{Test A}. The ROM dimension and the dimension of the hyperreduction basis are the same as for \textsf{Test A}; therefore, the performance on the test set of parameters (see right figure) is similar, with around $0.2\%$ relative error. The FOM flux waveform is compared with that obtained using the ROM in~\Cref{fig:alg2_ventricle_flux}. The ROM is able to accurately capture the behaviour at both the test samples $\gamma = 0.0004$ and $\gamma = 0.0057$. Furthermore, the speedup is similar as in \textsf{Test A}, i.e., $49$x.

\Cref{fig:alg2_ventricle_Phi_FOM_ROM} shows the solution snapshots of the transmembrane potential $\Phi$ (for both the FOM and the ROM) evaluated at the time steps $t^{k} \in \{20, 120, 200, 280, 400, 480, 520\}$ milliseconds for the parameter $\gamma = 0.00617$ taken from the test set. An excellent agreement can be seen between the transmembrane potential of the FOM and that of the ROM.
\section{Conclusions and future perspectives}%
\label{sec:conclusion}
 Without doubt, any particular disease progression in the heart does not occur in just a few cardiac cycles. In fact, certain abnormal conditions, e.g., high blood pressure or artery disease, cause an observable pathology after a long period of time, e.g., several months or even years. Therefore, if one is interested in modelling the disease progression from its early to severe stages, an extremely high computational load must be faced. In this context, ROMs have enormous potential to enable fast and reliable computer simulations of continuously beating heart over a long time period for variations of the input or parameter values, a key requirement for cardiac digital twins. Our work has introduced an adaptive algorithm to obtain ROMs. The adaptive algorithm is driven by an a posteriori error estimator which iteratively samples the best parameter to update the ROM basis. Furthermore, we also make use of a radial basis surrogate to determine the best parameter regions to focus on, in a problem-specific manner. The results on three examples have shown that the ROMs obtained from our adaptive algorithm generalize very well to unseen test cases. Moreover, they yield good approximation quality and capture the entire range of physical behaviour caused by input/parameter variations. For the scroll-wave phenomenon, we note that while the approximation quality of the ROM is sufficient, the acceleration in the simulation time is modest. This is mainly due to the comparatively larger ROM size for this example. This outcome is a limitation of linear ROM approaches applied to problems where the dynamics is dominated by convection or wave-like phenomena. To address this, our future work will investigate nonlinear projection-based approaches and machine learning approaches to obtain ROMs of smaller dimension. As an additional follow-up work, we plan to extend the introduced ROM methodology to excitation-contraction problems in the heart tissue where one often has large computational loads.



\addcontentsline{toc}{section}{References}
\bibliographystyle{myplainurl}
\bibliography{refs}

\begin{thebibliography}{10}

\bibitem{aliev+panfilov96}
R.~R. Aliev and A.~V. Panfilov.
\newblock A simple two-variable model of cardiac excitation.
\newblock {\em Chaos, Solitons \& Fractals}, 7:293--301, 1996.

\bibitem{AscRW95}
U.~M. Ascher, S.~J. Ruuth, and B.~T.~R. Wetton.
\newblock Implicit-explicit methods for time-dependent partial differential
  equations.
\newblock {\em SIAM Journal of Numerical Analysis}, 32(3):797--823, 1995.
\newblock \href {https://doi.org/10.1137/0732037} {\path{doi:10.1137/0732037}}.

\bibitem{Ballarin2016}
F.~Ballarin, E.~Faggiano, S.~Ippolito, A.~Manzoni, A.~Quarteroni, G.~Rozza, and
  R.~Scrofani.
\newblock Fast simulations of patient-specific haemodynamics of coronary artery
  bypass grafts based on a {POD}-{G}alerkin method and a vascular shape
  parametrization.
\newblock {\em Journal of Computational Physics}, 315:609--628, 2016.
\newblock \href {https://doi.org/10.1016/j.jcp.2016.03.065}
  {\path{doi:10.1016/j.jcp.2016.03.065}}.

\bibitem{morBarF22}
J.~Barnett and C.~Farhat.
\newblock Quadratic approximation manifold for mitigating the {K}olmogorov
  barrier in nonlinear projection-based model order reduction.
\newblock {\em Journal of Computational Physics}, 464:Paper No. 111348, 20,
  2022.
\newblock \href {https://doi.org/10.1016/j.jcp.2022.111348}
  {\path{doi:10.1016/j.jcp.2022.111348}}.

\bibitem{morBarFM23}
J.~Barnett, C.~Farhat, and Y.~Maday.
\newblock Neural-network-augmented projection-based model order reduction for
  mitigating the {K}olmogorov barrier to reducibility.
\newblock {\em Journal of Computational Physics}, 492:Paper No. 112420, 20,
  2023.
\newblock \href {https://doi.org/10.1016/j.jcp.2023.112420}
  {\path{doi:10.1016/j.jcp.2023.112420}}.

\bibitem{morBarMNetal04}
M.~Barrault, Y.~Maday, N.~C. Nguyen, and A.~T. Patera.
\newblock An `empirical interpolation' method: application to efficient
  reduced-basis discretization of partial differential equations.
\newblock {\em Comptes Rendus Mathematique}, 339:667--672, 2004.
\newblock \href {https://doi.org/10.1016/j.crma.2004.08.006}
  {\path{doi:10.1016/j.crma.2004.08.006}}.

\bibitem{Bonomi2017}
D.~Bonomi, A.~Manzoni, and A.~Quarteroni.
\newblock A matrix {DEIM} technique for model reduction of nonlinear
  parametrized problems in cardiac mechanics.
\newblock {\em Computer Methods in Applied Mechanics and Engineering},
  324:300--326, 2017.
\newblock \href {https://doi.org/10.1016/j.cma.2017.06.011}
  {\path{doi:10.1016/j.cma.2017.06.011}}.

\bibitem{morCagMS19}
N.~Cagniart, Y.~Maday, and B.~Stamm.
\newblock Model order reduction for problems with large convection effects.
\newblock In {\em Contributions to {P}artial {D}ifferential {E}quations and
  {A}pplications}, volume~47 of {\em Computer Methods in Applied Sciences},
  pages 131--150. Springer, Cham, 2019.

\bibitem{cansiz+dal+etal15}
B.~Cans{\i}z, H.~Dal, and M.~Kaliske.
\newblock An orthotropic viscoelastic material model for passive myocardium:
  Theory and algorithmic treatment.
\newblock {\em Computer Methods in Biomechanics and Biomedical Engineering},
  18:1160--1172, 2015.

\bibitem{cansiz+dal+etal17}
B.~Cans{\i}z, H.~Dal, and M.~Kaliske.
\newblock Computational cardiology: A modified hill model to describe the
  electro-visco-elasticity of the myocardium.
\newblock {\em Computer Methods in Applied Mechanics and Engineering},
  315:434--466, 2017.

\bibitem{cansiz+sveric+etal18}
B.~Cans{\i}z, K.~Sveric, K.~Ibrahim, R.~H. Strasser, A.~Linke, and M.~Kaliske.
\newblock Towards predictive computer simulations in cardiology: Finite element
  analysis of personalized heart models.
\newblock {\em ZAMM - Journal of Applied Mathematics and Mechanics /
  Zeitschrift für Angewandte Mathematik und Mechanik}, 98:2155--2176, 2018.

\bibitem{Cansiz2022}
B.~Cans{\i}z and M.~Kaliske.
\newblock A comparative study of fully implicit staggered and monolithic
  solution methods. {P}art {I}: {C}oupled bidomain equations of cardiac
  electrophysiology.
\newblock {\em Journal of Computational and Applied Mathematics}, 407:114021,
  21, 2022.
\newblock \href {https://doi.org/10.1016/j.cam.2021.114021}
  {\path{doi:10.1016/j.cam.2021.114021}}.

\bibitem{cansiz+woodworth+etal21}
B.~Cans{\i}z, L.~A. Woodworth, and M.~Kaliske.
\newblock A simple phenomenological approach for myocardial contraction:
  formulation, parameter sensitivity study and applications in organ level
  simulations.
\newblock {\em Mechanics of Soft Materials}, 3:1--28, 2021.

\bibitem{cansiz+kaliske22}
B.~Cansız and M.~Kaliske.
\newblock A comparative study of fully implicit staggered and monolithic
  solution methods. {P}art {I}: Coupled bidomain equations of cardiac
  electrophysiology.
\newblock {\em Journal of Computational Applied Mathamatics}, 407:114021, 2022.

\bibitem{Chak20}
N.~Chamakuri and P.~K\"{u}gler.
\newblock A coupled monodomain solver with optimal memory usage for the
  simulation of cardiac wave propagation.
\newblock {\em Applied Mathematics and Computation}, 378:125212, 15, 2020.
\newblock \href {https://doi.org/10.1016/j.amc.2020.125212}
  {\path{doi:10.1016/j.amc.2020.125212}}.

\bibitem{morChaS10}
S.~Chaturantabut and D.~C. Sorensen.
\newblock Nonlinear model reduction via discrete empirical interpolation.
\newblock {\em SIAM Journal on Scientific Computing}, 32:2737--2764, 2010.
\newblock \href {https://doi.org/10.1137/090766498}
  {\path{doi:10.1137/090766498}}.

\bibitem{morChe23}
S.~Chellappa.
\newblock {\em A Posteriori Error Estimation and Adaptivity for Model Order
  Reduction of Large-Scale Systems}.
\newblock {D}issertation, Otto-von-Guericke-Universit{\"a}t, Magdeburg,
  Germany, 2023.
\newblock \href {https://doi.org/http://dx.doi.org/10.25673/101396}
  {\path{doi:http://dx.doi.org/10.25673/101396}}.

\bibitem{morCheFB19a}
S.~Chellappa, L.~Feng, and P.~Benner.
\newblock Adaptive basis construction and improved error estimation for
  parametric nonlinear dynamical systems.
\newblock {\em International Journal for Numerical Methods in Engineering},
  121(23):5320--5349, 2020.
\newblock \href {https://doi.org/10.1002/nme.6462}
  {\path{doi:10.1002/nme.6462}}.

\bibitem{morCheFB22}
S.~Chellappa, L.~Feng, and P.~Benner.
\newblock An adaptive sampling approach for the reduced basis method.
\newblock In {\em Realization and Model Reduction of Dynamical Systems - A
  Festschrift in Honor of the 70th Birthday of {T}hanos {A}ntoulas}, pages
  137--155. Springer, Cham, 2022.
\newblock \href {https://doi.org/10.1007/978-3-030-95157-3_8}
  {\path{doi:10.1007/978-3-030-95157-3_8}}.

\bibitem{Franzone2014}
P.~Colli~Franzone, L.~F. Pavarino, and S.~Scacchi.
\newblock {\em Mathematical Cardiac Electrophysiology}, volume~13 of {\em
  MS\&A. Modeling, Simulation and Applications}.
\newblock Springer, Cham, 2014.
\newblock \href {https://doi.org/10.1007/978-3-319-04801-7}
  {\path{doi:10.1007/978-3-319-04801-7}}.

\bibitem{Corral-Acero_EHJ2020}
J.~Corral-Acero, F.~Margara, M.~Marciniak, C.~Rodero, F.~Loncaric, Y.~Feng,
  A.~Gilbert, J.~F. Fernandes, H.~A. Bukhari, A.~Wajdan, M.~V. Martinez, M.~S.
  Santos, M.~Shamohammdi, H.~Luo, P.~Westphal, P.~Leeson, P.~DiAchille,
  V.~Gurev, M.~Mayr, L.~Geris, P.~Pathmanathan, T.~Morrison, R.~Cornelussen,
  F.~Prinzen, T.~Delhaas, A.~Doltra, M.~Sitges, E.~J. Vigmond, E.~Zacur,
  V.~Grau, B.~Rodriguez, E.~W. Remme, S.~Niederer, P.~Mortier, K.~McLeod,
  M.~Potse, E.~Pueyo, A.~Bueno-Orovio, and P.~Lamata.
\newblock {The ‘Digital Twin’ to enable the vision of precision
  cardiology}.
\newblock {\em European Heart Journal}, 41(48):4556--4564, 2020.
\newblock \href {https://doi.org/10.1093/eurheartj/ehaa159}
  {\path{doi:10.1093/eurheartj/ehaa159}}.

\bibitem{einthoven12}
W.~Einthoven.
\newblock The different forms of the human electrocardiogram and their
  signification.
\newblock {\em The Lancet}, 179:853--861, 1912.

\bibitem{morFenCB23}
L.~Feng, S.~Chellappa, and P.~Benner.
\newblock A posteriori error estimation for model order reduction of parametric
  systems.
\newblock preprint, 2023.
\newblock \href {https://doi.org/10.21203/rs.3.rs-3410762/v1}
  {\path{doi:10.21203/rs.3.rs-3410762/v1}}.

\bibitem{fitzhugh61}
R.~Fitzhugh.
\newblock Impulses and physiological states in theoretical models of nerve
  membrane.
\newblock {\em Biophysical Journal}, 1:445--466, 1961.

\bibitem{Fresca2020}
S.~Fresca, A.~Manzoni, L.~Ded{\`e}, and A.~Quarteroni.
\newblock Deep learning-based reduced order models in cardiac
  electrophysiology.
\newblock {\em {PLOS} {ONE}}, 15(10):e0239416, 2020.

\bibitem{Fresca2021}
S.~Fresca, A.~Manzoni, L.~Ded{\`e}, and A.~Quarteroni.
\newblock {POD}-enhanced deep learning-based reduced order models for the
  real-time simulation of cardiac electrophysiology in the left atrium.
\newblock {\em Frontiers in Physiology}, 12:679076, 2021.

\bibitem{Gerbeau2015}
J.-F. Gerbeau, D.~Lombardi, and E.~Schenone.
\newblock Reduced order model in cardiac electrophysiology with approximated
  {L}ax pairs.
\newblock {\em Advances in Computational Mathematics}, 41(5):1103--1130, 2015.
\newblock \href {https://doi.org/10.1007/s10444-014-9393-9}
  {\path{doi:10.1007/s10444-014-9393-9}}.

\bibitem{goktepe+kuhl09}
S.~G\"oktepe and E.~Kuhl.
\newblock Computational modeling of cardiac electrophysiology: a novel finite
  element approach.
\newblock {\em International Journal for Numerical Methods in Engineering},
  79:156--178, 2009.

\bibitem{goktepe+wong+etal10}
S.~G{\"o}ktepe, J.~Wong, and E.~Kuhl.
\newblock Atrial and ventricular fibrillation: computational simulation of
  spiral waves in cardiac tissue.
\newblock {\em Archive of Applied Mechanics}, 80:569--580, 2010.

\bibitem{GreU19}
C.~Greif and K.~Urban.
\newblock Decay of the {K}olmogorov {$N$}-width for wave problems.
\newblock {\em Applied Mathematics Letters}, 96:216--222, 2019.
\newblock \href {https://doi.org/10.1016/j.aml.2019.05.013}
  {\path{doi:10.1016/j.aml.2019.05.013}}.

\bibitem{morGre05}
M.~Grepl.
\newblock {\em Reduced-basis approximation a posteriori error estimation for
  parabolic partial differential equations}.
\newblock PhD thesis, Massachussetts Institute of Technology (MIT), Cambridge,
  USA, 2005.
\newblock URL: \url{http://dspace.mit.edu/handle/1721.1/7582}.

\bibitem{morHaaO08}
B.~Haasdonk and M.~Ohlberger.
\newblock Reduced basis method for finite volume approximations of parametrized
  linear evolution equations.
\newblock {\em {ESAIM}: Mathematical Modelling and Numerical Analysis},
  42:277--302, 2008.
\newblock \href {https://doi.org/10.1051/m2an:2008001}
  {\path{doi:10.1051/m2an:2008001}}.

\bibitem{morHaaO11}
B.~Haasdonk and M.~Ohlberger.
\newblock Efficient reduced models and a posteriori error estimation for
  parametrized dynamical systems by offline/online decomposition.
\newblock {\em Mathematical and Computer Modelling of Dynamical Systems},
  17(2):145--161, 2011.
\newblock \href {https://doi.org/10.1080/13873954.2010.514703}
  {\path{doi:10.1080/13873954.2010.514703}}.

\bibitem{morHesRS16}
J.~S. Hesthaven, G.~Rozza, and B.~Stamm.
\newblock {\em Certified Reduced Basis Methods for Parametrized Partial
  Differential Equations}.
\newblock SpringerBriefs in Mathematics. Springer International Publishing,
  2016.
\newblock \href {https://doi.org/10.1007/978-3-319-22470-1}
  {\path{doi:10.1007/978-3-319-22470-1}}.

\bibitem{hodgkin+huxley90}
A.~Hodgkin and A.~Huxley.
\newblock A quantitative description of membrane current and its application to
  conduction and excitation in nerve.
\newblock {\em Bulletin of Mathematical Biology}, 52:25--71, 1990.

\bibitem{keating+sanguinetti01}
M.~T. Keating and M.~C. Sanguinetti.
\newblock Molecular and cellular mechanisms of cardiac arrhythmias.
\newblock {\em Cell}, 104:569--580, 2001.

\bibitem{Khan_2022}
R.~Khan and K.~T. Ng.
\newblock Numerical study of {POD}-{G}alerkin-{DEIM} reduced order modeling of
  cardiac monodomain formulation.
\newblock {\em Biomedical Physics \& Engineering Express}, 8:015012, 2021.
\newblock \href {https://doi.org/10.1088/2057-1976/ac3c0b}
  {\path{doi:10.1088/2057-1976/ac3c0b}}.

\bibitem{kotikanyadanam+goktepe+etal10}
M.~Kotikanyadanam, S.~Göktepe, and E.~Kuhl.
\newblock Computational modeling of electrocardiograms: A finite element
  approach toward cardiac excitation.
\newblock {\em International Journal for Numerical Methods in Biomedical
  Engineering}, 26:524--533, 2010.

\bibitem{Krishnamoorthi2013}
S.~Krishnamoorthi, M.~Sarkar, and W.~S. Klug.
\newblock Numerical quadrature and operator splitting in finite element methods
  for cardiac electrophysiology.
\newblock {\em Internatial Journal of Numerical Methods in Biomedical
  Engineering}, 29:1243--1266, 2013.
\newblock \href {https://doi.org/10.1002/cnm.2573}
  {\path{doi:10.1002/cnm.2573}}.

\bibitem{Loewe2022}
A.~Loewe, P.~Mart{\'i}nez~D{\'i}az, C.~Nagel, and J.~S{\'a}nchez.
\newblock Cardiac {D}igital {T}win {M}odeling.
\newblock In T.~Jadczyk, G.~Caluori, A.~Loewe, and K.~S. Golba, editors, {\em
  Innovative Treatment Strategies for Clinical Electrophysiology}, pages
  111--134. Springer Nature, 2022.
\newblock \href {https://doi.org/10.1007/978-981-19-6649-1_7}
  {\path{doi:10.1007/978-981-19-6649-1_7}}.

\bibitem{Manzoni2018}
A.~Manzoni, D.~Bonomi, and A.~Quarteroni.
\newblock {\em Reduced Order Modeling for Cardiac Electrophysiology and
  Mechanics: New Methodologies, Challenges and Perspectives}, pages 115--166.
\newblock Springer, 2018.
\newblock \href {https://doi.org/10.1007/978-3-319-96649-6_6}
  {\path{doi:10.1007/978-3-319-96649-6_6}}.

\bibitem{nagumo+arimoto+etal62}
J.~Nagumo, S.~Arimoto, and S.~Yoshizawa.
\newblock An active pulse transmission line simulating nerve axon.
\newblock {\em Proceedings of the {IRE}}, 50:2061--2070, 1962.

\bibitem{niederer_computational_2019}
S.~A. Niederer, J.~Lumens, and N.~A. Trayanova.
\newblock Computational models in cardiology.
\newblock {\em Nature Reviews Cardiology}, 16:100--111, 2019.
\newblock \href {https://doi.org/10.1038/s41569-018-0104-y}
  {\path{doi:10.1038/s41569-018-0104-y}}.

\bibitem{Pagani2018}
S.~Pagani, A.~Manzoni, and A.~Quarteroni.
\newblock Numerical approximation of parametrized problems in cardiac
  electrophysiology by a local reduced basis method.
\newblock {\em Computer Methods in Applied Mechanics and Engineering},
  340:530--558, 2018.
\newblock \href {https://doi.org/10.1016/j.cma.2018.06.003}
  {\path{doi:10.1016/j.cma.2018.06.003}}.

\bibitem{Pathmanathan2011}
P.~Pathmanathan, G.~R. Mirams, J.~Southern, and J.~P. Whiteley.
\newblock The significant effect of the choice of ionic current integration
  method in cardiac electro-physiological simulations.
\newblock {\em International Journal for Numerical Methods in Biomedical
  Engineering}, 27:1751--1770, 2011.
\newblock \href {https://doi.org/10.1002/cnm.1438}
  {\path{doi:10.1002/cnm.1438}}.

\bibitem{Peirlincketal2021}
M.~Peirlinck, F.~Costabal, J.~Yao, J.~Guccione, S.~Tripathy, Y.~Wang,
  D.~Ozturk, P.~Segars, T.~Morrison, S.~Levine, and E.~Kuhl.
\newblock Precision medicine in human heart modeling: Perspectives, challenges,
  and opportunities.
\newblock {\em Biomechanics and Modeling in Mechanobiology (online)},
  20:803--831, 2021.
\newblock \href {https://doi.org/10.1007/s10237-021-01421-z}
  {\path{doi:10.1007/s10237-021-01421-z}}.

\bibitem{Pfaller2020}
M.~R. Pfaller, M.~Cruz~Varona, J.~Lang, C.~Bertoglio, and W.~A. Wall.
\newblock Using parametric model order reduction for inverse analysis of large
  nonlinear cardiac simulations.
\newblock {\em International Journal for Numerical Methods in Biomedical
  Engineering}, 36:e3320, 27, 2020.
\newblock \href {https://doi.org/10.1002/cnm.3320}
  {\path{doi:10.1002/cnm.3320}}.

\bibitem{Puwal-Roth2007}
S.~Puwal and B.~J. Roth.
\newblock Forward {E}uler stability of the bidomain model of cardiac tissue.
\newblock {\em IEEE Transactions on Biomedical Engineering}, 54:951--953, 2007.
\newblock \href {https://doi.org/10.1109/TBME.2006.889204}
  {\path{doi:10.1109/TBME.2006.889204}}.

\bibitem{Zhilin1999_operatorsplit1}
Z.~Qu and A.~Garfinkel.
\newblock An advanced algorithm for solving partial differential equation in
  cardiac conduction.
\newblock {\em IEEE Transactions on Biomedical Engineering}, 46(9):1166--1168,
  1999.
\newblock \href {https://doi.org/10.1109/10.784149}
  {\path{doi:10.1109/10.784149}}.

\bibitem{morQuaMN16}
A.~Quarteroni, A.~Manzoni, and F.~Negri.
\newblock {\em {R}educed {B}asis {M}ethods for {P}artial {D}ifferential
  {E}quations}, volume~92 of {\em La Matematica per il 3+2}.
\newblock Springer International Publishing, 2016.
\newblock \href {https://doi.org/10.1007/978-3-319-15431-2}
  {\path{doi:10.1007/978-3-319-15431-2}}.

\bibitem{Quarteroni_AN2017}
A.~Quarteroni, A.~Manzoni, and C.~Vergara.
\newblock The cardiovascular system: mathematical modelling, numerical
  algorithms and clinical applications.
\newblock {\em Acta Numerica}, 26:365--590, 2017.
\newblock \href {https://doi.org/10.1017/S0962492917000046}
  {\path{doi:10.1017/S0962492917000046}}.

\bibitem{QuaDM19}
A.~Quarteroni, L.~Dede', A.~Manzoni, and C.~Vergara.
\newblock {\em Mathematical Modelling of the Human Cardiovascular System: Data,
  Numerical Approximation, Clinical Applications}.
\newblock Cambridge Monographs on Applied and Computational Mathematics.
  Cambridge University Press, 2019.
\newblock \href {https://doi.org/10.1017/9781108616096}
  {\path{doi:10.1017/9781108616096}}.

\bibitem{rogers02}
J.~M. Rogers.
\newblock Wave front fragmentation due to ventricular geometry in a model of
  the rabbit heart.
\newblock {\em Chaos {(Woodbury,} {N.Y.)}}, 12:779---87, 2002.

\bibitem{rogers+mcculloch94a}
J.~M. Rogers and A.~D. McCulloch.
\newblock A {collocation-Galerkin} finite element model of cardiac action
  potential propagation.
\newblock {\em {IEEE} Transactions on Biomedical Engineering}, 41:743--757,
  1994.

\bibitem{Sachse_2005}
F.~B. Sachse.
\newblock {\em Computational Cardiology: Modeling of Anatomy,
  Electrophysiology, and Mechanics}.
\newblock Lecture Notes in Computer Science. Springer, 2005.
\newblock \href {https://doi.org/10.1007/10.1007/b96841}
  {\path{doi:10.1007/10.1007/b96841}}.

\bibitem{sundnes+lines+etal05}
J.~Sundnes, G.~T. Lines, and A.~Tveito.
\newblock An operator splitting method for solving the bidomain equations
  coupled to a volume conductor model for the torso.
\newblock {\em Mathematical Biosciences}, 194:233--248, 2005.

\bibitem{taylor_feap}
R.~L. Taylor.
\newblock {FEAP} - {F}inite {E}lement {A}nalysis {P}rogram, 2020.
\newblock URL: \url{http://www.ce.berkeley/feap}.

\bibitem{tentusscher+noble+etal04}
K.~H.~W. {ten Tusscher}, D.~Noble, P.~J. Noble, and A.~V. Panfilov.
\newblock A model for human ventricular tissue.
\newblock {\em American Journal of Physiology-Heart and Circulatory
  Physiology}, 286:1573--1589, 2004.

\bibitem{wong+goktepe+etal11}
J.~Wong, S.~G\"oktepe, and E.~Kuhl.
\newblock Computational modeling of electrochemical coupling: A novel finite
  element approach towards ionic models for cardiac electrophysiology.
\newblock {\em Computer Methods in Applied Mechanics and Engineering},
  200:3139--3158, 2011.

\bibitem{Woodworth21}
L.~A. Woodworth, B.~Cans{\i}z, and M.~Kaliske.
\newblock A numerical study on the effects of spatial and temporal
  discretization in cardiac electrophysiology.
\newblock {\em International Journal for Numerical Methods in Biomedical
  Engineering}, 37(5):e3443, 23, 2021.
\newblock \href {https://doi.org/10.1002/cnm.3443}
  {\path{doi:10.1002/cnm.3443}}.

\bibitem{woodworth+cansiz+etal22}
L.~A. Woodworth, B.~Cans{\i}z, and M.~Kaliske.
\newblock Balancing conduction velocity error in cardiac electrophysiology
  using a modified quadrature approach.
\newblock {\em International Journal for Numerical Methods in Biomedical
  Engineering}, 38(5):e3589, 2022.

\bibitem{YangVeneziani2017}
H.~Yang and A.~Veneziani.
\newblock Efficient estimation of cardiac conductivities via {POD}-{DEIM} model
  order reduction.
\newblock {\em Appl. Numer. Math.}, 115:180--199, 2017.
\newblock \href {https://doi.org/10.1016/j.apnum.2017.01.006}
  {\path{doi:10.1016/j.apnum.2017.01.006}}.

\bibitem{Yeetal22}
D.~Ye, P.~Zun, V.~Krzhizhanovskaya, and A.~G. Hoekstra.
\newblock Uncertainty quantification of a three-dimensional in-stent restenosis
  model with surrogate modelling.
\newblock {\em Journal of the Royal Society Interface}, 19(187):20210864, 2022.
\newblock \href {https://doi.org/10.1098/rsif.2021.0864}
  {\path{doi:10.1098/rsif.2021.0864}}.

\bibitem{morZhaFLetal15}
Y.~Zhang, L.~Feng, S.~Li, and P.~Benner.
\newblock An efficient output error estimation for model order reduction of
  parametrized evolution equations.
\newblock {\em {SIAM} Journal on Scientific Computing}, 37(6):B910--B936, 2015.
\newblock \href {https://doi.org/10.1137/140998603}
  {\path{doi:10.1137/140998603}}.

\end{thebibliography}
  
\end{document}